\documentclass[11pt]{article}
\usepackage{amsmath,amsfonts,amssymb}
\newtheorem{theorem}{Theorem}

\newtheorem{condition}[theorem]{Condition}

\newtheorem{corollary}[theorem]{Corollary}

\newtheorem{definition}[theorem]{Definition}

\newtheorem{lemma}[theorem]{Lemma}

\newtheorem{proposition}[theorem]{Proposition}
\newtheorem{remark}[theorem]{Remark}

\newenvironment{proof}[1][Proof]{\textbf{#1.} }{\ \rule{0.5em}{0.5em}}

\newcommand{\R}{\mathbb{R}}
\newcommand{\RR}{\mathbb{R}}
\newcommand{\eps}{\varepsilon}
\newcommand{\expect}{E}
\newcommand{\prob}{P}
\newcommand{\const}{\textrm{const}}

\newcommand{\cov}{\textrm{Cov}}

\newcommand{\wZ}{ Z^-}

\begin{document}

\title{On the regularity of 
stochastic currents,\\ fractional Brownian motion\\ and applications 
to a turbulence model}
\author{
\\
  { Franco Flandoli}
           \\
  {\small\it Dipartimento di Matematica Applicata,
Universit\`{a} di Pisa}    \\[-0.1cm]
  {\small\it Via Bonanno 25B, I-56126 Pisa, Italia}          \\[-0.1cm]
 {\small {\tt flandoli@dma.unipi.it}    
}   
\and
  \\[-0.1cm]
{  Massimiliano Gubinelli } \\
 {\small\it    Laboratoire de Math\'ematiques,  
 Universit\'e de Paris-Sud }    \\[-0.1cm]
  {\small\it  B\^atiment 425,   F-91405 Orsay Cedex, France}      \\[-0.1cm]
 {\small  {\tt massimiliano.gubinelli@math.u-psud.fr}   
}   
 \and
 \\[-0.1cm]
   { Francesco Russo }              \\
   {\small\it Institut Galil\'ee, Math\'ematiques, Universit\'e Paris
     13 } 
   \\[-0.1cm]
   {\small\it 99, av. J.-B. Cl\'ement - F-93430 Villetaneuse, France} 
         \\[-0.1cm]
   {\small  {\tt russo@math.univ-paris13.fr}}   
}
\date{March  2007}
\maketitle
\begin{abstract}
We study the pathwise regularity of the map
$$
\varphi \mapsto I(\varphi) = \int_0^T \langle \varphi(X_t), dX_t \rangle
$$  
where $\varphi$ is a vector function on $\R^d$ belonging to some Banach
space $V$, $X$ is a stochastic process and the integral is
some version of a stochastic integral defined via regularization. A \emph{stochastic current} is
a continuous version of this map, seen as a random element of the
topological dual of $V$. We give  sufficient conditions 
for the current to live in some Sobolev space of distributions and we 
provide elements to conjecture that those are also necessary.
 Next we verify the sufficient conditions when the process $X$ is a $d$-dimensional
fractional Brownian motion (fBm);
we identify regularity in Sobolev spaces for fBm with Hurst
index $H \in (1/4,1)$. Next we provide some results about general
Sobolev regularity of Brownian currents. Finally we discuss
applications to a  model of random vortex filaments in
turbulent fluids.


\end{abstract}



\section{Introduction}

We consider stochastic integrals, loosely speaking of the form
\begin{equation} \label{ECurr}
I\left(  \varphi\right)  =\int_{0}^{T}\left\langle \varphi\left(
X_{t}\right)  ,dX_{t}\right\rangle, 
\end{equation}
where $\left(  X_{t}\right)  $ is a Wiener process or a fractional Brownian
motion with Hurst parameter $H$ in a certain range. We are interested in the
pathwise continuity properties with respect to $\varphi$: we would like to
establish that the random generalized field $\varphi\mapsto I\left(  \varphi\right)  $ has
a version that is a.s. continuous in $\varphi$ in certain topologies. In the
language of geometric measure theory, such a property means that the
stochastic integral defines pathwise a \textit{current}, with the regularity
specified by the topologies that we have found.

This problem is motivated by the study of  fluidodynamical models. In \cite{Ffil},
in the study of the energy of a vortex filament naturally appear 
some stochastic double integral related to Wiener process
\begin{equation} \label{EDint}
 \int_{[0,T]^2} f(X_s - X_t) dX_s dX_t,
\end{equation}
where $f(x) =  K_{\alpha}\left(  x\right)$
where $K_{\alpha}\left(  x \right)  $ is the kernel of the pseudo-differential operator $(1-\Delta)^{-\alpha}$
(precise definitions will be given in section
\ref{preliminaries}). $f$ is therefore a continuous singular function at zero.

The difficulty there comes from the appearance of anticipating integrands
and from the singularity of $f$ at zero.
 \cite{Ffil} gives sense to this integral in some Stratonovich sense.
Moreover that paper explores the connection with self-intersection local time considered
for instance by J.-F. Le Gall in \cite{gall}.

The work~\cite{Nua} considers a similar double integral in the case of fractional Brownian 
motion with Hurst index $H  > \frac{1}{2}$ using Malliavin-Skorohod anticipating calculus.

A natural approach is to interpret previous double integral as a symmetric (or eventually) forward integral
in the framework of  stochastic calculus via regularization, see \cite{RVSem} for a 
survey. 
We recall that when $X$ is a semimartingale,  forward  (respectively symmetric) integral
$\int_0^t \varphi(X) d^- X$ (resp. $\int_0^t \varphi(X) d^\circ X$)   coincides with the corresponding It\^o
(resp. Stratonovich)  integral.
The double stochastic integral considered by \cite{Ffil} coincides in
fact with the symmetric integral introduced here. 
So (\ref{EDint}) can be interpreted as 
\begin{equation} \label{EDintsym}
 \int_{[0,T]^2} f(X_s - X_t) d^\circ X_s d^\circ X_t,
\end{equation}

In this paper, $X$ will be a fractional Brownian motion
with Hurst index $ H > 1/4$ 
   but a complete study of the existence of integrals (\ref{EDintsym})
will be not yet  performed here because of heavy technicalities.
We will only essentially  consider 
their regularized versions. Now, those double integrals are naturally
in correspondence with 
currents related to  $I$  defined in (\ref{ECurr}). The investigation
of those currents is strictly related to ``pathwise stochastic
calculus'' in the spirit of \emph{rough paths} theory by T. Lyons and
coauthors,  \cite{Lyons, qianlyons}.
Here we aim at exploring the ``pathwise character'' of 
stochastic integrals via regularization.
 A first step in this direction was done
in \cite{GradNo}, where the authors showed that forward integrals
 of the type
 $\int_0^T \varphi(X) d^- X$, when $X$ is a one-dimensional semimartingale 
or a fractional Brownian motion with Hurst index $H > \frac{1}{2}$,
can be regarded as a.s. uniform approximations of their regularization $I^-_\eps(\varphi)$ (see section 3),
instead of the usual  convergence in probability.

This analysis of currents related to stochastic integrals was 
started in \cite{FGGT} using an approach based on 
spectral analysis. The approach presented here is not based on Fourier
transform and contains new general ideas with respect to \cite{FGGT}.
 Informally speaking, it is based on the formula
\begin{equation}
\int_{0}^{T}\left\langle \varphi\left(  X_{t}\right)  ,dX_{t}\right\rangle
=\int_{\mathbb{R}^{d}}\left\langle (1-\Delta)^{\alpha}\varphi\left(  x\right)
,\int_{0}^{T}K_{\alpha}\left(  x-X_{t}\right)  dX_{t}\right\rangle
dx .\label{decoupling}%
\end{equation}
This
formula decouples $\varphi$ and $X$ and replaces the problem of the pathwise
dependence of $I\left(  \varphi\right)  $ on the infinite dimensional
parameter $\varphi$ with the problem of the pathwise dependence of $\int
_{0}^{T}K_{\alpha}\left(  x-X_{t}\right)  dX_{t}$ on the finite dimensional
parameter $x\in\mathbb{R}^{d}$. Another form of  decoupling is also
one of the ingredients of the Fourier approach of \cite{FGGT} but the
novelty here is that we can take better advantages from the properties
of the underlying process (like, for example, the existence of a
density). Moreover formula~(\ref{decoupling}) produces at least two
new results.

First, we can treat in an essentially optimal way the case of fractional
Brownian motion, making use of its Gaussian properties. For $H>1/2$ results in
this direction can be extracted from the estimates proved in \cite{Nua} again
by spectral analysis. However, with the present approach we may treat the case
$H\in(1/4,1/2)$ as well.

Second, in the case of the Brownian motion, we may work with functions
$\varphi$ in the Sobolev spaces of Banach type $H_{p}^{\alpha}$, with $p>1$,
instead of only the Hilbert topologies $H_{2}^{\alpha}$ considered in~
\cite{FGGT}, with the great advantage that it is sufficient to ask less
differentiability on $\varphi$ (any $\alpha>1$ suffices), at the price of a
larger $p$ (depending on $\alpha$ and the space dimension). In this way we may
cover, for instance, the class $\varphi\in C^{1,\varepsilon}$ treated
in~\cite{Lyons}, see also \cite{Gubi}; the approach here is entirely different
and does not rely on rough paths, see Remark \ref{Rrough}.

Finally, we apply these ideas to random vortex filaments. In the case of the
fractional Brownian motion we prove new results about the finiteness of the
kinetic energy of the filaments (such a property is expected to be
linked to the
regularity of the pathwise current). In the case of the Brownian motion
optimal conditions for a finite energy were already proved in~\cite{Ffil,FGub}%
, while a sufficient condition when $H>1/2$ has been found in~ \cite{Nua}. The
results of the present work provide new regularity properties of the random
filaments, especially for the parameter range $H \in(1/4,1/2) $.

\section{Generalities}\label{sec:hilbert}

\subsection{Stochastic currents}

Let $\left(  X_{t}\right)  $ be a stochastic process such that $X_0 = 0$ a.s.
 on a probability space
$\left(  \Omega,\mathcal{F},P\right)  $ with values in
 $\mathbb{R}^{d}$. Let $T > 0$. 
Let
$V $ be a Banach space of vector fields $\varphi:\mathbb{R}^{d}\rightarrow
\mathbb{R}^{d}$ and $\mathcal{D}\subset V$ be a dense subset. Assume that a
stochastic integral
\[
I\left(  \varphi\right)  =\int_{0}^{T}\left\langle \varphi\left(
X_{t}\right)  ,dX_{t}\right\rangle
\]
is well defined, in a suitable sense (It\^{o}, etc.), for every $\varphi
\in\mathcal{D}$. Our first aim is to define it for every $\varphi\in V$. In
addition, we would like to prove that it has a pathwise redefinition according
to the following:

\begin{definition}
\label{def:pathwise} The family of r.v. $\{I(\varphi)\}_{\varphi\in V}$ has a
{\bf pathwise redefinition}  on $V$ if there exists a measurable mapping $\xi
:\Omega\rightarrow V^{\prime}$ such that for every $\varphi\in\mathcal{D}$
\begin{equation}
I\left(  \varphi\right)  \left(  \omega\right)  =\left(  \xi\left(
\omega\right)  \right)  \left(  \varphi\right)  \text{ for }P\text{-a.e.
}\omega\in\Omega.\label{pathwise}%
\end{equation}
\end{definition}

Then, if we succeed in our objective,

\begin{itemize}
\item  for every $\varphi\in V$ we consider the r.v. $\omega\mapsto\left(
\xi\left(  \omega\right)  \right)  \left(  \varphi\right)  $ as the definition
of the stochastic integral $I\left(  \varphi\right)  $ (now extended to the
class $\varphi\in V$ )

\item  for $P$-a.e. $\omega\in\Omega$, we consider the linear continuous
mapping $\varphi\mapsto\left(  \xi\left(  \omega\right)  \right)  \left(
\varphi\right)  $ as a pathwise redefinition of stochastic integral on $V$.
\end{itemize}

Formally, the candidate for $\xi$ is the expression
\[
\xi\left(  x\right)  =\int_{0}^{T}\delta\left(  x-X_{t}\right)  dX_{t}.
\]
where $\delta$ is here the $d$-dimensional Dirac measure.
Indeed, always formally,
\begin{align*}
\xi\left(  \varphi\right)   &  =\int_{\mathbb{R}^{d}}\langle\xi\left(
x\right) , \varphi\left(  x\right) \rangle dx=\int_{0}^{T}\left\langle \left(
\int_{\mathbb{R}^{d}}\delta\left(  x-X_{t}\right)  \varphi\left(  x\right)
dx\right)  , dX_{t}\right\rangle \\
&  =\int_{0}^{T}\langle\varphi( X_{t}), dX_{t}\rangle=I(\varphi).
\end{align*}
We remark that this  viewpoint is inspired by the theory of
currents; with other methods (spectral ones) it was developed in~\cite{FGGT}.

\subsection{Decoupling by duality}

As we said in the introduction, our approach is based on a proper rigorous
version of formula (\ref{decoupling}). One way to interpret it by the
following duality argument, that we describe only at a formal level.

Let $W$ be another Banach space and $\Lambda:V\rightarrow W$ be an
isomorphism. Proceeding formally as above we have
\begin{align*}
\int_{0}^{T}\left\langle \varphi\left(  X_{t}\right)  ,dX_{t}\right\rangle
_{\mathbb{R}^{d}}  &  =\left\langle \varphi,\xi\right\rangle _{V,V^{\prime}%
}=\left\langle \Lambda^{-1}\Lambda\varphi,\xi\right\rangle _{V,V^{\prime}}\\
&  =\left\langle \Lambda\varphi,\left(  \Lambda^{-1}\right)  ^{\ast}%
\xi\right\rangle _{W,W^{\prime}}%
\end{align*}
(notice that $\Lambda^{-1}:W\rightarrow V$, $\left(  \Lambda^{-1}\right)
^{\ast}:V^{\prime}\rightarrow W^{\prime}$).
Our aim
essentially amounts to prove that
$
\left(  \Lambda^{-1}\right)  ^{\ast}\xi:\Omega\rightarrow W^{\prime}%
$
is a well defined random variable.

This reformulation becomes useful if the spaces $W$, $W^{\prime}$ are easier
to handle than $V$, $V^{\prime}$, and the operator $\left(  \Lambda
^{-1}\right)  ^{\ast}$ has a kernel $K(x,y)$ as an operator in function
spaces:
\[
\left(  \left(  \Lambda^{-1}\right)  ^{\ast}f\right)  \left(  x\right)  =\int
K\left(  x,y\right)  f(y)dy.
\]
In such a case, formally
\[%
\begin{split}
\left(  \left(  \Lambda^{-1}\right)  ^{\ast}\xi\right)  \left(  x\right)   &
=\int K\left(  x,y\right)  \left(  \int_{0}^{T}\delta\left(  y-X_{t}\right)
dX_{t}\right)  dy\\
&  =\int_{0}^{T}\left(  \int K\left(  x,y\right)  \delta\left(  y-X_{t}%
\right)  dy\right)  dX_{t}\\
&  =\int_{0}^{T}K\left(  x,X_{t}\right)  dX_{t}.
\end{split}
\]
Below we make a rigorous version of this representation by choosing $V=H_{p}%
^{\alpha}\left(  \mathbb{R}^{d}\right)  $, $W=L^{p}\left(  \mathbb{R}%
^{d}\right)  $, $\Lambda=(1-\Delta)^{\frac{\alpha}{2}}$, $K\left(  x,y\right)
=$ $K_{\alpha/2}\left(  x-y\right)  $ (notations are given in the next section).

\subsection{Rigorous setting}
\label{preliminaries}
Denote by
$S\left(  \mathbb{R}^{d}\right)  $ the space of rapidly decreasing infinitely
differentiable \textit{vector fields} $\varphi:\mathbb{R}^{d}\rightarrow
\mathbb{R}^{d}$, by $S^{\prime}\left(  \mathbb{R}^{d}\right)  $ its dual (the
space of tempered distributional fields) and by $\mathcal{F}$ the Fourier
transform
\[
\left(  \mathcal{F}\varphi\right)  \left(  \ell\right)  =\left(  2\pi\right)
^{-d/2}\int_{\mathbb{R}^{d}}e^{-i\left\langle x,\ell\right\rangle }%
\varphi(x)dx,\quad\ell\in\mathbb{R}^{d}%
\]
which is an isomorphism in both $S\left(  \mathbb{R}^{d}\right)  $ and
$S^{\prime}\left(  \mathbb{R}^{d}\right)  $. Let $\mathcal{F}^{-1}$ denote the
inverse Fourier transform. For every $s\in\mathbb{R}$, let $\Lambda
_{s}:S^{\prime}\left(  \mathbb{R}^{d}\right)  \rightarrow S^{\prime}\left(
\mathbb{R}^{d}\right)  $ be the pseudo-differential operator defined as
\[
\Lambda_{s}\varphi=\mathcal{F}^{-1}\left(  1+\left|  \ell\right|  ^{2}\right)
^{\frac{s}{2}}\mathcal{F}\varphi.
\]
We shall also denote it by $(1-\Delta)^{\frac{s}{2}}$.

Let $H_{p}^{s}\left(  \mathbb{R}^{d}\right)  $, with $p>1$ and $s\in\mathbb{R}
$, be the Sobolev space of vector fields $\varphi\in S^{\prime}\left(
\mathbb{R}^{d}\right)  $ such that
\[
\left\|  \varphi\right\|  _{H_{p}^{s}}^{p}:=\int_{\mathbb{R}^{d}}%
|(1-\Delta)^{\frac{s}{2}}\varphi(x)|^{p}dx<\infty,
\]
see~\cite{Triebel}, sec. 2.3.3,  where the definition
chosen here for brevity is given as a characterization. From the very
definitions of $\Lambda_{s}$ and $H_{p}^{s}\left(  \mathbb{R}%
^{d}\right)  $, the operator $\Lambda_{s}$ is an isomorphism from $H_{p}%
^{s}\left(  \mathbb{R}^{d}\right)  $ onto $L^{p}\left(  \mathbb{R}^{d}\right)
$ (the Lebesgue space of $p$-integrable \textit{vector fields}).

Another fact often used in the paper is that the dual space $\left(  H_{p}%
^{s}\left(  \mathbb{R}^{d}\right)  \right)  ^{\prime}$ is $H_{p^{\prime}}%
^{-s}\left(  \mathbb{R}^{d}\right)  $:
\[
\left(  H_{p}^{s}\left(  \mathbb{R}^{d}\right)  \right)  ^{\prime
}=H_{p^{\prime}}^{-s}\left(  \mathbb{R}^{d}\right)  ,\quad\frac{1}{p}+\frac
{1}{p^{\prime}}=1,
\]
see \cite{Triebel}, section 2.6.1. Moreover, being $\Lambda_{s}$ an
isomorphism from $H_{p}^{s}\left(  \mathbb{R}^{d}\right)  $ onto $L^{p}\left(
\mathbb{R}^{d}\right)  $, its dual operator $\Lambda_{s}^{\star}$ is an
isomorphism from $L^{p^{\prime}}\left(  \mathbb{R}^{d}\right)  $ onto
$H_{p^{\prime}}^{-s}\left(  \mathbb{R}^{d}\right)  $.

It is known that negative fractional powers of a positive selfadjoint operator
$A$ in a Hilbert space $H$, such that $-A$ generates the semigroup $T(t)$,
have the representation
\[
A^{-\alpha}=\frac{1}{\Gamma\left(  \alpha\right)  }\int_{0}^{\infty}%
t^{\alpha-1}T\left(  t\right)  dt,
\]
where $\Gamma$ is the standard Gamma function, see \cite{Pa}, formula (6.9). Taking $A=$ $(1-\Delta)$ in the Hilbert space
$H=L^{2}\left(  \mathbb{R}^{d}\right)  $, we have
\[
\left(  T\left(  t\right)  \varphi\right)  \left(  x\right)  =\left(  4\pi
t\right)  ^{-d/2}\int_{\mathbb{R}^{d}}e^{-\frac{\left|  x-y\right|  ^{2}}%
{4t}-t}\varphi(y)dy
\]
and thus, for $\alpha>0$,
\begin{align*}
\left(  (1-\Delta)^{-\alpha}\varphi\right)  (x) &  =\frac{\left(  4\pi
t\right)  ^{-d/2}}{\Gamma\left(  \alpha\right)  }\int_{0}^{\infty}t^{\alpha
-1}\int_{\mathbb{R}^{d}}e^{-\frac{\left|  x-y\right|  ^{2}}{4t}-t}%
\varphi(y)dydt\\
&  =\int_{\mathbb{R}^{d}}\left[  \frac{1}{\Gamma\left(  \alpha\right)  \left(
4\pi\right)  ^{d/2}}\int_{0}^{\infty}t^{\alpha-1-\frac{d}{2}}e^{-\frac{\left|
x-y\right|  ^{2}}{4t}-t}dt\right]  \varphi(y)dy.
\end{align*}
In fact this formula can be proved more elementarily from the definition of
$(1-\Delta)^{-\alpha}\varphi$ and the formula
\[
\lambda^{-\alpha}=\frac{1}{\Gamma\left(  \alpha\right)  }\int_{0}^{\infty
}t^{\alpha-1}e^{-\lambda t}dt,
\]
then taking $\lambda=1+\left|  \ell\right|  ^{2}$ and the Fourier transform of
the Gaussian density. This fact implies that the operator $(1-\Delta
)^{-\alpha}$, which originally is an isomorphism between $H_{2}^{-2\alpha
}\left(  \mathbb{R}^{d}\right)  $ and $L^{2}\left(  \mathbb{R}^{d}\right)  $,
considered by restriction as a bounded linear operator in $L^{2}\left(
\mathbb{R}^{d}\right)  $, has a kernel $K_{\alpha}\left(  .\right)  $,
\begin{equation}
\left(  (1-\Delta)^{-\alpha}\varphi\right)  (x)=\int_{\mathbb{R}^{d}}%
K_{\alpha}\left(  x-y\right)  \varphi(y)dy\label{eq kernel}%
\end{equation}
given by
\[
K_{\alpha}\left(  x\right)  =\frac{1}{\Gamma\left(  \alpha\right)  \left(
4\pi\right)  ^{d/2}}\int_{0}^{\infty}t^{\alpha-\frac{d}{2}}e^{-\frac{\left|
x\right|  ^{2}}{4t}-t}\frac{dt}{t}.
\]
The following estimates are not optimized as far as the exponential decay is
concerned; we just state a version sufficient for our purposes.
The proof of the two lemmas before are in Appendix A.

\begin{lemma}
\label{lemma su K}There exists positive constants
$c_{\alpha,d}$, $C_{\alpha,d}$ such that:
\begin{itemize}
\item[1)] 
For $0<\alpha<\frac{d}{2}$, we have
\begin{equation} \label{eq:asympK}
K_{\alpha}\left(  x\right)  =\left|  x\right|  ^{2\alpha-d}\rho\left(
x\right)
\end{equation}
where
$
c_{\alpha,d}e^{-2\left|  x\right|  ^{2}}\leq\rho\left(  x\right)  \leq
C_{\alpha,d}e^{-\frac{\left|  x\right|  }{8}}%
$;
\item[2)] 
 For $\alpha
>\frac{d}{2}$, we have
\begin{equation} \label{eq:asympK2new}
c_{\alpha,d} e^{-\frac{\left|  x\right|  ^{2}}{4}}\leq K_{\alpha
}\left(  x\right)  \leq C_{\alpha,d} e^{-\frac{\left|  x\right|  }{8}}%
\end{equation}
for two positive constants $c_{\alpha,d}^{\prime}$, $C_{\alpha,d}^{\prime}$.
Moreover
\begin{equation}
  \label{eq:asympK2}
K_{\alpha}(x) = K_\alpha(0) - \rho'(x) |x|^{-d+2\alpha} \ge 0   
\end{equation}
with $0 < K_\alpha(0) < \infty$ and where
$
c_{\alpha,d}e^{-2\left|  x\right|  ^{2}}\leq\rho'\left(  x\right)  \leq
C_{\alpha,d}e^{-\frac{\left|  x\right|  }{8}}%
$;

\item[3)]
Finally, when $\alpha=d/2$ we have
\begin{equation}
  \label{eq:log-asymp-K}
 K_\alpha(x) \le C_{\alpha,d} \ \log |x| e^{- a_\alpha  \ |x|} 
\end{equation}
where $a_\alpha$ is another positive constant.
\end{itemize}
\end{lemma}
\begin{remark} \label{RsuK}
In particular $\rho$ and $\rho'$ are bounded.
\end{remark}

In the applications we will need also some control on $\Delta K_\alpha(x)$ which is provided by the next lemma.

\begin{lemma}
\label{lemma:deltaK}
It holds that $-\Delta K_{\alpha}(x) = K_\alpha(x) - K_{\alpha-1}(x)$. Then, when $\alpha<d/2+1$,
\begin{equation}
 |-\Delta K_{\alpha}(x)| \le \rho^{\prime\prime}(x)|x|^{2\alpha-d-2}%
\label{eq:asympDeltaK}%
\end{equation}
where $\rho^{\prime\prime}$ is positive, bounded above, locally bounded away from
zero below and depends on $\alpha$.
\end{lemma}

\subsection{Regularity of stochastic currents}

With these notations and preliminaries in mind, we may state a first rigorous
variant of formula (\ref{decoupling}). Given a continuous stochastic process
$\left(  X_{t}\right)  _{t\geq0}$ on $\left(  \Omega,\mathcal{F},P\right)  $
with values in $\mathbb{R}^{d}$, given $\varepsilon>0 $, let $\left(
D_{\varepsilon}X_{t}\right)  _{t\geq0}$ be any one of the following discrete
derivatives:
\[
\frac{X_{t+\varepsilon}-X_{t}}{\varepsilon},\quad\frac{X_{t}-X_{t-\varepsilon
}}{\varepsilon},\quad\frac{X_{t+\varepsilon}-X_{t-\varepsilon}}{\varepsilon}%
\]
where we understand that $X_{t-\varepsilon}=0$ for $t<\varepsilon$. The
following integral
\[
I_{\varepsilon}\left(  \varphi\right)  =\int_{0}^{T}\left\langle
\varphi\left(  X_{t}\right)  ,D_{\varepsilon}X_{t}\right\rangle dt
\]
is well defined $P$-a.s. as a classical Riemann integral, at least for every
continuous vector field $\varphi$.

\begin{lemma} \label{L5}
Given  $\alpha,\varepsilon>0$,
with probability one the function $t\mapsto K_{\alpha/2}\left(  x-X_{t}%
\right)  D_{\varepsilon}X_{t}$ is integrable for a.e. $x\in\mathbb{R}^{d}$,
the function
\[
\eta_{\varepsilon}\left(  x\right)  :=\int_{0}^{T}K_{\alpha/2}\left(
x-X_{t}\right)  D_{\varepsilon}X_{t}dt
\]
is in $L^{1}\left(  \mathbb{R}^{d}\right)  $ and 
for any $\varphi\in S\left(  \mathbb{R}^{d}\right)  $ we have
\begin{equation}
\int_{0}^{T}\left\langle \varphi\left(  X_{t}\right)  ,D_{\varepsilon}%
X_{t}\right\rangle dt=\int_{\mathbb{R}^{d}}\left\langle (1-\Delta)^{\alpha
/2}\varphi\left(  x\right)  ,\int_{0}^{T}K_{\alpha/2}\left(  x-X_{t}\right)
D_{\varepsilon}X_{t}dt\right\rangle dx.\label{decoupling epsilon}%
\end{equation}
\end{lemma}

\begin{proof}
\bigskip Notice that $(1-\Delta)^{\alpha/2}\varphi\in S\left(  \mathbb{R}%
^{d}\right)  $ and, by (\ref{eq kernel}),
\begin{align*}
\varphi=(1-\Delta)^{-\alpha/2}(1-\Delta)^{\alpha/2}\varphi
=\int_{\mathbb{R}^{d}}K_{\alpha/2}\left(  \cdot-x\right)  \left[
(1-\Delta)^{\alpha/2}\varphi\right]  (x)dx.
\end{align*}
Thus
\[
\int_{0}^{T}\left\langle \varphi\left(  X_{t}\right)  ,D_{\varepsilon}%
X_{t}\right\rangle dt=\int_{0}^{T}\left\langle \int_{\mathbb{R}^{d}}%
K_{\alpha/2}\left(  X_{t}-x\right)  \left[  (1-\Delta)^{\alpha/2}%
\varphi\right]  (x)dx,D_{\varepsilon}X_{t}\right\rangle dt.
\]
Denote by $\widehat{K}_{\alpha/2}\left(  x\right)  $ the function equal to
$K_{\alpha/2}\left(  x\right)  $ for $x\neq0$, infinite for $x=0$. 
Suppose for a moment that
\begin{equation}
P\left(  \int_{0}^{T}\int_{\mathbb{R}^{d}}\widehat{K}_{\alpha/2}\left(
X_{t}-x\right)  dxdt<\infty\right)  =1. \label{eq integrability}%
\end{equation}
Then the integrability properties stated in the lemma
will hold.
Since
\[
  \sup_{  \left(  t,x\right)  \in\left[  0,T\right]
  \times\mathbb{R}^{d}} \left \vert
\left\langle \left[  (1-\Delta)^{\alpha/2}\varphi\right]
(x), D_{\varepsilon}X_{t}\right\rangle \right \vert <\infty \quad {a.s.},
\]
using   Fubini theorem  we will get (\ref{decoupling epsilon})
 (notice that $K_{\alpha/2}\left(
X_{t}-x\right)  =K_{\alpha/2}\left(  x-X_{t}\right)  $).

Thus we have only to
prove (\ref{eq integrability}). Since $K_{\alpha/2}$ is positive, we may apply again
Fubini theorem and analyze
$ \int_{0}^{T}(  \int_{\mathbb{R}^{d}}\widehat{K}_{\alpha/2}(
X_{t}-x)  dx)  dt
$.
 But we have, for every $y\in\mathbb{R}^{d}$%
\[
\int_{\mathbb{R}^{d}}\widehat{K}_{\alpha/2}\left(  y-x\right)  dx = \int_{\mathbb{R}^{d}}\widehat{K}_{\alpha/2}\left( x\right) dx.
\]
This quantity   is finite for every $\alpha>0$, from the estimates of lemma \ref{lemma
su K}. The proof is now complete.
\end{proof}

\bigskip
Below we need a criterion to decide when $\eta_{\varepsilon}$ (defined in the
previous lemma) belongs to $L^{2}\left(  \mathbb{R}^{d}\right)  $. It is thus
useful to introduce the following condition which ensures the existence
of the representation given in Lemma~\ref{ltech5} below.

\begin{condition}
\label{cond:B}
$
\int_{0}^{T}\int_{0}^{T} \widehat 
K_{\alpha}\left(  X_{t}-X_{s}\right)  dtds<\infty
$
$\prob$-a.s.
\end{condition}

\begin{remark}
 $\int_{0}^{T}%
\int_{0}^{T} \widehat K_{\alpha}\left(  X_{t}-X_{s}\right)  dtds<\infty$, 
implies 
that the function $\left(  t,s\right)  \mapsto X_{t}-X_{s}$ is different from
zero except possibly on a zero measure set of $\left[  0,T\right]  ^{2}$, and that the
well-defined function $\left(  t,s\right)  \mapsto \widehat K_{\alpha}\left(
X_{t}-X_{s}\right)  $ is Lebesgue integrable on $\left[  0,T\right]  ^{2}$.
From now on the notation $  \widehat 
K_{\alpha} $ will simply be replaced by 
$K_{\alpha}$.
\end{remark}

\begin{lemma} \label{ltech5}
Under Condition~\ref{cond:B} we have the following double integral
representation for the norm of $\eta_\eps$ defined in Lemma \ref{L5}:
\[
\left\|  \eta_{\varepsilon}\right\|  _{L^{2}(\mathbb{R}^{d})}^{2}=\int_{0}%
^{T}\int_{0}^{T}K_{\alpha}\left(  X_{t}-X_{s}\right)  \left\langle
D_{\varepsilon}X_{t},D_{\varepsilon}X_{s}\right\rangle dtds.
\]
\end{lemma}
\begin{proof}
We have
\begin{align*}
\left\|  \eta_{\varepsilon}\right\|  _{L^{2}(\mathbb{R}^{d})}^{2}  &
=\int_{\mathbb{R}^{d}}\left\langle \eta_{\varepsilon}\left(  x\right)
,\eta_{\varepsilon}\left(  x\right)  \right\rangle dx\\
& =\int_{\mathbb{R}^{d}}\left\langle \int_{0}^{T}K_{\alpha/2}\left(
x-X_{t}\right)  D_{\varepsilon}X_{t}dt,\int_{0}^{T}K_{\alpha/2}\left(
x-X_{s}\right)  D_{\varepsilon}X_{s}ds\right\rangle dx\\
& \int_{0}^{T}\int_{0}^{T}\int_{\mathbb{R}^{d}}K_{\alpha/2}\left(
x-X_{t}\right)  K_{\alpha/2}\left(  x-X_{s}\right)  \left\langle
D_{\varepsilon}X_{t},D_{\varepsilon}X_{s}\right\rangle dxdtds
\end{align*}
if we can apply Fubini Theorem. Then it is sufficient to use the property
\[
\int_{\mathbb{R}^{d}}K_{\alpha/2}\left(  x-y\right)  K_{\alpha/2}\left(
x-z\right)  dx=K_{\alpha}\left(  y-z\right)  .
\]
Since the process $X$ is continuous, both $(t,s) \mapsto D_\eps X_s D_\eps X_t$ is a continuous  two-parameter process which on $[0,T]^2$ is a.s. bounded. Then a sufficient condition to apply Fubini Theorem is
\[
\int_{0}^{T}\int_{0}^{T}\int_{\mathbb{R}^{d}}K_{\alpha/2}\left(
x-X_{t}\right)  K_{\alpha/2}\left(  x-X_{s}\right)  dxdtds<\infty.
\]
Since the integrand is positive, it is equal to
\begin{align*}
& \int_{0}^{T}\int_{0}^{T}\left(  \int_{\mathbb{R}^{d}}K_{\alpha/2}\left(
x-X_{t}\right)  K_{\alpha/2}\left(  x-X_{s}\right)  dx\right)  dtds
=\int_{0}^{T}\int_{0}^{T}K_{\alpha}\left(  X_{t}-X_{s}\right)  dtds.
\end{align*}

Invoking condition~\ref{cond:B} we can conclude the proof.
\end{proof}

\bigskip
The double integral representation of the norm of $\eta_\eps$ will play a major r\^ole in the following, so we introduce the notation
$$
Z_{\alpha,\eps} := \int_{0}%
^{T}\int_{0}^{T}K_{\alpha}\left(  X_{t}-X_{s}\right)  \left\langle
D_{\varepsilon}X_{t},D_{\varepsilon}X_{s}\right\rangle dtds.
$$

\begin{lemma}
\label{lemma:ordering}
Assume Condition~\ref{cond:B} holds for any $\alpha \ge \overline{\alpha}$. Then the function $\alpha \mapsto Z_{\alpha,\eps} \ge 0$ is decreasing for $\alpha \ge \overline{\alpha}$.  
\end{lemma}
\begin{proof}
Denote $\eta_{\alpha,\eps}(x) = \int_0^T K_{\alpha/2}(x-X_t) D_\eps X_t  dt $ making explicit the dependence on $\alpha$. It is not difficult to prove that, if $\overline{\alpha} \le \alpha \le \beta $ we have
$
\eta_{\beta,\eps} = (1-\Delta)^{(\alpha-\beta)/2} \eta_{\alpha,\eps}
$.
Then
$Z_{\beta,\eps} = \|\eta_{\beta,\eps}\|^2_{L^2(\RR^d)} = \|(1-\Delta)^{(\alpha-\beta)/2}\eta_{\alpha,\eps}\|^2_{L^2(\RR^d)}   \le  \|\eta_{\alpha,\eps}\|^2_{L^2(\RR^d)}= Z_{\alpha,\eps}$, since being $\alpha-\beta\le 0$ , the operator $(1-\Delta)^{(\alpha-\beta)/2}$ has a norm bounded by one.
\end{proof}

\bigskip

We will assume below the following condition on the convergence of the regularized integrals.
\begin{condition}
\label{cond:A} For every $\varphi\in C_{0}^{\infty}\left(  \mathbb{R}%
^{d}\right)  $, $I_{\varepsilon}\left(  \varphi\right)  $ converges in
probability to some r.v., denoted by $I\left(  \varphi\right)  $.
\end{condition}

Under Condition~\ref{cond:A}, the mapping $\varphi\mapsto I\left(
\varphi\right)  $ is a priori defined only on $C_{0}^{\infty}\left(
\mathbb{R}^{d}\right)  $ with values in the set $L^{0}\left(  \Omega\right)  $
of random variables. Its extension to $\varphi\in H_{2}^{\alpha}\left(
\mathbb{R}^{d}\right)  $ is a result of the next theorem.

\begin{theorem}
\label{th:l2result} Assume Conditions~\ref{cond:B} and~\ref{cond:A},
and the a priori bound
\[
\sup_{\varepsilon\in\left(  0,1\right)  }E\left[  \int_{0}^{T}\int_{0}%
^{T}K_{\alpha}\left(  X_{t}-X_{s}\right)  \left\langle D_{\varepsilon}%
X_{t},D_{\varepsilon}X_{s}\right\rangle dtds\right]  <\infty.
\]
Then:

i) the mappings $\varphi\in C_{0}^{\infty}\left(  \mathbb{R}^{d}\right)
\mapsto I_{\varepsilon}\left(  \varphi\right)  ,I\left(  \varphi\right)  \in
L^{0}\left(  \Omega\right)  $ take values in $L^{2}\left(  \Omega\right)  $
and extend (uniquely) to linear continuous mappings from $H_{2}^{\alpha
}\left(  \mathbb{R}^{d}\right)  $ to $L^{2}\left(  \Omega\right)  $. Moreover,
for every $\varphi\in H_{2}^{\alpha}\left(  \mathbb{R}^{d}\right)  $,
$I_{\varepsilon}\left(  \varphi\right)  \rightarrow I\left(  \varphi\right)  $
in probability and in $L^{2-\delta}\left(  \Omega\right)  $ for every
$\delta>0$.

ii) In addition, there exist random elements $\xi_{\varepsilon},\xi
:\Omega\rightarrow H_{2}^{-\alpha}\left(  \mathbb{R}^{d}\right)  $ (in fact
 belonging to $L^{2}(  \Omega;H_{2}^{-\alpha}(  \mathbb{R}^{d}))  $)
that constitute pathwise redefinitions of $I_{\varepsilon}$ and $I$ over the
functions $\varphi\in H_{2}^{\alpha}\left(  \mathbb{R}^{d}\right)  $, in the
sense of Definition~\ref{pathwise}.
\end{theorem}

\begin{proof}
\textbf{Step 1} (mean square results). By the assumptions and the previous
lemma we have
$
\sup_{\varepsilon\in\left(  0,1\right)  }E\left[  \left\|  \eta_{\varepsilon
}\right\|  _{L^{2}(\mathbb{R}^{d})}^{2}\right]  <\infty
$.
From (\ref{decoupling epsilon}), for $\varphi\in C_{0}^{\infty}\left(
\mathbb{R}^{d}\right)  $ we have
\[
I_{\varepsilon}\left(  \varphi\right)  =\left\langle \left(  1-\Delta\right)
^{\alpha/2}\varphi,\eta_{\varepsilon}\right\rangle _{L^{2}(\mathbb{R}^{d})}.
\]
Therefore, always for $\varphi\in C_{0}^{\infty}\left(  \mathbb{R}^{d}\right)
$,
\[
\left|  I_{\varepsilon}\left(  \varphi\right)  \right|  \leq\left\|
\eta_{\varepsilon}\right\|  _{L^{2}(\mathbb{R}^{d})}\left\|  \left(
1-\Delta\right)  ^{\alpha/2}\varphi\right\|  _{L^{2}\left(  \mathbb{R}%
^{d}\right)  }\leq C_{\alpha}\left\|  \eta_{\varepsilon}\right\|
_{L^{2}(\mathbb{R}^{d})}\left\|  \varphi\right\|  _{H_{2}^{\alpha}\left(
\mathbb{R}^{d}\right)  }%
\]%
\[
E\left[  \left|  I_{\varepsilon}\left(  \varphi\right)  \right|  ^{2}\right]
\leq C_{\alpha}\left\|  \varphi\right\|  _{H_{2}^{\alpha}\left(
\mathbb{R}^{d}\right)  }^{2}E\left[  \left\|  \eta_{\varepsilon}\right\|
_{L^{2}\left(  \mathbb{R}^{d}\right)  }^{2}\right]  \leq C_{\alpha}^{\prime
}\left\|  \varphi\right\|  _{H_{2}^{\alpha}\left(  \mathbb{R}^{d}\right)
}^{2}.
\]
\ 

Immediately we have $I_{\varepsilon}\left(  \varphi\right)  \in L^{2}\left(
\Omega\right)  $ for every $\varphi\in C_{0}^{\infty}\left(  \mathbb{R}%
^{d}\right)  $, and the mapping $\varphi\mapsto I_{\varepsilon}\left(
\varphi\right)  $ extends (uniquely by density) to a linear continuous mapping
from $H_{2}^{\alpha}\left(  \mathbb{R}^{d}\right)  $ to $L^{2}\left(
\Omega\right)  $.

Given $\varphi\in C_{0}^{\infty}(  \mathbb{R}^{d})  $, since
$I_{\varepsilon}(  \varphi)  \rightarrow I(  \varphi) $
in probability, uniform integrability arguments and $E[  |
I_{\varepsilon}(  \varphi)  |  ^{2}]  \leq C_{\alpha
}^{\prime}\|  \varphi\|  _{H_{2}^{\alpha}(  \mathbb{R}%
^{d})  }^{2}$, yield $I_{\varepsilon}(  \varphi)  \rightarrow
I(  \varphi)  $ in $L^{2-\delta}(  \Omega)  $ for every
$\delta>0$; moreover, it is not difficult to deduce $I(  \varphi)
\in L^{2}(  \Omega)  $ and
$
E[ |  I(  \varphi)  |  ^{2}]  \leq
C_{\alpha}^{\prime}\|  \varphi\|  _{H_{2}^{\alpha}(
\mathbb{R}^{d})  }^{2}
$.
As before, this implies that the mapping $\varphi\mapsto I\left(
\varphi\right)  $ extends uniquely to a linear continuous mapping from
$H_{2}^{\alpha}\left(  \mathbb{R}^{d}\right)  $ to $L^{2}\left(
\Omega\right)  $. Now, with these extensions, it is not difficult to show that
$I_{\varepsilon}\left(  \varphi\right)  \rightarrow I\left(  \varphi\right)  $
in $L^{2-\delta}\left(  \Omega\right)  $ for every $\delta>0$ also for every
for $\varphi\in H_{2}^{\alpha}\left(  \mathbb{R}^{d}\right)  $.

\medskip
\textbf{Step 2} (pathwise results). We still have to construct $\xi
_{\varepsilon}$ and $\xi$. Recalling that $\eta_{\varepsilon}\in L^{2}\left(
\Omega;L^{2}(\mathbb{R}^{d})\right)  $, $\xi_{\varepsilon}$ is simply defined
as $[  (  1-\Delta)  ^{\alpha/2}]  ^{\star}%
\eta_{\varepsilon}$, element of $L^{2}\left(  \Omega;H_{2}^{-\alpha}\left(
\mathbb{R}^{d}\right)  \right)  $, where $A^\star$ denotes the dual of an
operator $A$.

 To this end, recall that $\left(
1-\Delta\right)  ^{\alpha/2}$ is an isomorphism between $H_{2}^{\alpha}\left(
\mathbb{R}^{d}\right)  $ and $L^{2}(\mathbb{R}^{d})$, and thus the dual
operator $[ (  1-\Delta)  ^{\alpha/2}]  ^{\star}$ is
an isomorphism between the dual spaces $L^{2}(\mathbb{R}^{d})$ and
$H_{2}^{-\alpha}\left(  \mathbb{R}^{d}\right)  $ (we identify $L^{2}%
(\mathbb{R}^{d})$ with its dual).

The family $\left\{  \eta_{\varepsilon}\right\}  $ is bounded in $L^{2}\left(
\Omega;L^{2}(\mathbb{R}^{d})\right)  $, hence there exist a sequence
$\eta_{\varepsilon_{n}}$ weakly convergent to some $\eta$ in $L^{2}\left(
\Omega;L^{2}(\mathbb{R}^{d})\right)  $:
$
E\left\langle \eta_{\varepsilon_{n}},Y\right\rangle _{L^{2}\left(
\mathbb{R}^{d}\right)  }\rightarrow E\left\langle \eta,Y\right\rangle
_{L^{2}(\mathbb{R}^{d})}%
$
for every $Y\in L^{2}(  \Omega;L^{2}(\mathbb{R}^{d}))  $. We set
$\xi:=[  (  1-\Delta)  ^{\alpha/2}]  ^{\star}\eta$,
random element of $H_{2}^{-\alpha}\left(  \mathbb{R}^{d}\right)  $. We shall
see that this definition does not depend on the sequence $\varepsilon_{n}$. We
have to prove that for every $\varphi\in C_{0}^{\infty}\left(  \mathbb{R}%
^{d}\right)  $,
$
I\left(  \varphi\right)  \left(  \omega\right)  =\left(  \xi\left(
\omega\right)  \right)  \left(  \varphi\right)  \text{ for }P\text{-a.s.
}\omega\in\Omega.
$
Equivalently we have to prove that for every $\varphi\in C_{0}^{\infty}\left(
\mathbb{R}^{d}\right)  $

\begin{equation}
I\left(  \varphi\right)  \left(  \omega\right)  =
\eta
\left (\omega \right ) \left ( \left(
1-\Delta \right)  ^{\alpha/2} \varphi\right)   \text{ for
}P\text{-a.s. }\omega\in\Omega.\label{tobeproved}%
\end{equation}

We already know that $I_{\varepsilon}(\varphi)  =\langle
(  1-\Delta)  ^{\alpha/2}\varphi,\eta_{\varepsilon}\rangle
_{L^{2}(\mathbb{R}^{d})}$ for every $\varepsilon>0$. Choose $Y$ above of the
form $Y=F\left(  1-\Delta\right)  ^{\alpha/2}\varphi$ with generic $F\in
L^{2}\left(  \Omega\right)  $. Given $\varphi\in C_{0}^{\infty}\left(
\mathbb{R}^{d}\right)  $, we know that
\[
E\left[  F\left\langle \eta_{\varepsilon_{n}},\left(  1-\Delta\right)
^{\alpha/2}\varphi\right\rangle _{L^{2}(\mathbb{R}^{d})}\right]  \rightarrow
E\left[  F\left\langle \eta,\left(  1-\Delta\right)  ^{\alpha/2}%
\varphi\right\rangle _{L^{2}(\mathbb{R}^{d})}\right]
\]
for every $F\in L^{2}\left(  \Omega\right)  $. Hence
\[
E\left[  F I_{\varepsilon}\left(  \varphi\right)  \right]  \rightarrow E\left[
F\left\langle \eta,\left(  1-\Delta\right)  ^{\alpha/2}\varphi\right\rangle
_{L^{2}(\mathbb{R}^{d})}\right]
\]
but we also know that $I_{\varepsilon}\left(  \varphi\right)  \rightarrow
I\left(  \varphi\right)  $ in $L^{2-\delta}\left(  \Omega\right)  $ for every
$\delta>0$. We get
\[
E\left[  F I\left(  \varphi\right)  \right]  =E\left[  F \left\langle
\eta,\left(  1-\Delta\right)  ^{\alpha/2}\varphi\right\rangle _{L^{2}%
(\mathbb{R}^{d})}\right]
\]
at least for every bounded random variable $F$, hence (\ref{tobeproved}) holds
true. This also implies that the definition of $\xi$ does not depend on the
sequence $\varepsilon_{n}$. The proof is complete.
\end{proof}
\bigskip

To state a possible converse of Th.~\ref{th:l2result} it is useful to introduce a weaker version of
Def.~\ref{def:pathwise}.

\begin{definition}
\label{def:non-pathwise}  Let $A\in\mathcal{F}$ be such that $P(A)>0$. We say that $\{I(\varphi)\}_{\varphi\in V}$ has a
{\bf pathwise redefinition on} $A$ if there exists a measurable mapping $\xi:A\rightarrow
V^{\prime}$ such that for every $\varphi\in C_{0}^{\infty}\left(
\mathbb{R}^{d}\right)  $ equation~(\ref{pathwise}) holds true for $P$-a.e.
$\omega\in A$.
\end{definition}

\begin{theorem}
\label{th:neg-res-l2} Assume Conditions~\ref{cond:B} and~\ref{cond:A},
and the a priori bound%
\[
E\left[  \int_{0}^{T}\int_{0}^{T}K_{\alpha}\left(  X_{t}-X_{s}\right)
\left\langle D_{\varepsilon}X_{t},D_{\varepsilon}X_{s}\right\rangle
dtds\right]  <\infty
\]
for every $\varepsilon>0$. If there exists $A\in\mathcal{F}$ with $P(A)>0$
such that $\{I(\varphi)\}_{\varphi\in V}$ has a pathwise redefinition on $A$, then
\[
\underset{\varepsilon\rightarrow0}{\lim\sup}\int_{0}^{T}\int_{0}^{T}K_{\alpha
}\left(  X_{t}-X_{s}\right)  \left\langle D_{\varepsilon}X_{t},D_{\varepsilon
}X_{s}\right\rangle dtds<\infty\quad\text{$P$-a.s. $\omega\in A$}.
\]
\end{theorem}

\begin{proof}
The first part of step 1 of the previous proof is still valid (except for the
uniformity in $\varepsilon$ of the constants). Thus in particular $E[
\|  \eta_{\varepsilon}\|  _{L^{2}(\mathbb{R}^{d})}^{2}]
<\infty$, $I_{\varepsilon}(  \varphi)  =\langle (
1-\Delta)  ^{\alpha/2}\varphi,\eta_{\varepsilon}\rangle $,
$|  I_{\varepsilon}(  \varphi) |  \leq C_{\alpha
,\varepsilon}\|  \varphi\|  _{H_{2}^{\alpha}(  \mathbb{R}%
^{d})  }^{2}$, $I_{\varepsilon}(  \varphi)  \in L^{2}(
\Omega)  $ for every $\varphi\in C_{0}^{\infty}(  \mathbb{R}%
^{d})  $ and the mapping $\varphi\mapsto I_{\varepsilon}\left(
\varphi\right)  $ extends to a linear continuous mapping from $H_{2}^{\alpha
}\left(  \mathbb{R}^{d}\right)  $ to $L^{2}\left(  \Omega\right)  $.

We know, by Condition~\ref{cond:A} and the existence of $\xi$, that for every
$\varphi\in C_{0}^{\infty}\left(  \mathbb{R}^{d}\right)  $%
\[
\left\langle \left(  1-\Delta\right)  ^{\alpha/2}\varphi,\eta_{\varepsilon
}\left(  \omega\right)  \right\rangle _{L^{2}\left(  \Omega\right)  }%
\underset{\varepsilon\rightarrow0}{\rightarrow}\left(  \xi\left(
\omega\right)  \right)  \left(  \varphi\right)  \text{ for }P\text{-a.e.
}\omega\in A.
\]
One can find a countable set $ {\cal D}  \subset C_{0}^{\infty}\left(  {\mathbb R}^{d}\right)  $
 with the following two properties: i) $\left(  1-\Delta\right)
^{\alpha/2}{\cal D}$ is dense in $L^{2}(\mathbb{R}^{d})$ and ii) for $P$-a.e.
$\omega\in A$
\begin{equation}
\left\langle \left(  1-\Delta\right)  ^{\alpha/2}\varphi,\eta_{\varepsilon
}\left(  \omega\right)  \right\rangle _{L^{2}\left(  \Omega\right)  }%
\underset{\varepsilon\rightarrow0}{\rightarrow}\left(  \xi\left(
\omega\right)  \right)  \left(  \varphi\right)  \text{ for every }\varphi\in
\cal{D}.\label{limite aus}%
\end{equation}

Let us prove the claim by contradiction. Assume there is $A^{\prime}%
\in\mathcal{F}$, $A^{\prime}\subset A$, $P(A')>0$,
 such that
\[
\underset{\varepsilon\rightarrow0}{\lim\sup}\int_{0}^{T}\int_{0}^{T}K_{\alpha
}\left(  X_{t}-X_{s}\right)  \left\langle D_{\varepsilon}X_{t},D_{\varepsilon
}X_{s}\right\rangle dtds=\infty\quad\text{for every }\omega\in A^{\prime
}\text{.}%
\]
By  lemma \ref{ltech5},
$
\underset{\varepsilon\rightarrow0}{\lim\sup}\left\|  \eta_{\varepsilon}\left(
\omega\right)  \right\|  _{L^{2}(\mathbb{R}^{d})}^{2}=\infty$
for
every $\omega\in A^{\prime}$.

Consequently, there is a subset   $A^{\prime\prime}$ of   $A^{\prime}$
 with $P(A^{\prime\prime} ) > 0$
such that (\ref{limite aus}) holds true for every $\omega\in
A^{\prime\prime}$. Thus, given $\omega\in A^{\prime\prime}$, there is an
infinitesimal sequence $\left\{  \varepsilon_{n}\left(  \omega\right)
\right\}  $ such that
$
\lim_{n\rightarrow\infty}\left\|  \eta_{\varepsilon_{n}\left(  \omega\right)
}\left(  \omega\right)  \right\|  _{L^{2}(\mathbb{R}^{d})}^{2}=\infty
$
but at the same time
\[
\lim_{n\rightarrow\infty}\left\langle \left(  1-\Delta\right)  ^{\alpha
/2}\varphi,\eta_{\varepsilon_{n}\left(  \omega\right)  }\left(  \omega\right)
\right\rangle _{L^{2}\left(  \Omega\right)  }=\left(  \xi\left(
\omega\right)  \right)  \left(  \varphi\right)
\]
for every $\varphi\in \cal{D}$. This is impossible because of the density of
$\left(  1-\Delta\right)  ^{\alpha/2} \cal{D} $ in $L^{2}(\mathbb{R}^{d})$. The proof
is complete.
\end{proof}

\begin{definition}
We say that a.s. a stochastic current $I$ does not belong to $V'$ if there is no
 $A \in \mathcal{F}$, $\prob(A) > 0$ such that $\{ I(\varphi) \}_{\varphi \in V}$ has a pathwise redefinition on $A$.   
\end{definition}

\begin{corollary}
\label{cor:neg-res}
Under the conditions of Theorem ~\ref{th:neg-res-l2}, if  
\[
\underset{\varepsilon\rightarrow0}{\lim\sup}\int_{0}^{T}\int_{0}^{T}K_{\alpha
}\left(  X_{t}-X_{s}\right)  \left\langle D_{\varepsilon}X_{t},D_{\varepsilon
}X_{s}\right\rangle dtds= +\infty\quad\text{$P$-a.s.}
\]
the current $I$ does not belong to $H^\alpha_2(\RR^d)$. 
\end{corollary}

\section{Application to the fractional Brownian motion}

We recall here that a fractional Brownian motion (fBm) $B = (B_t)$ with 
Hurst index $H\in\left(  0,1\right )$, is a Gaussian mean-zero real process whose 
covariance function is given by 
$$ Cov(B_s,B_t)  = \frac{1}{2}(|s|^{2H}+|t|^{2H}-|s-t|^{2H}),\,(s,t)\in \mathbb{R}_+^{2}.$$
This process has been widely studied: for some recent developments, we point to 
\cite{em} as a relevant monograph.
For instance we recall that when $H =  \frac{1}{2}$, $B$ is a classical Wiener process.
Its trajectories are H\"older continuous with respect to any parameter $\gamma < H$.
Recall that 
$$ E( \vert B_t - B_s \vert^2) = \vert t - s \vert ^{2H}.$$ 
For our stochastic integral redefinition, we have chosen the framework of
stochastic calculus via regularization; for survey about the topic and
recent developments, see \cite{RVSem}.

Let $X=\left(  X^{1},...,X^{d}\right)  $ be a $d$-dimensional fractional Brownian motion with
Hurst index $H\in\left(  0,1\right)$,
i.e. an $\mathbb{R}^d$-valued process whose components  $X^{1},...,X^{d}$ are real independent fractional
Brownian motions. In this section we study the
regularity of  the current generated by $X$ using the results of the
previous section and in particular Theorem~\ref{th:l2result} which
gives sufficient conditions for regularity in the Hilbert spaces
$H^{\alpha}_{2}(\RR^d)$.

We will consider the \emph{symmetric} and \emph{forward} integrals,
respectively defined as the limit in probability, as $\eps \to 0$, of
$$
I^\circ_\eps(\varphi) := \int_0^T < \varphi(X_t), \frac{X_{t+\eps}-X_{t-\eps}}{2\eps}> dt
$$
and
$$
 I^-_\eps(\varphi) := \int_0^T <\varphi(X_t), \frac{X_{t+\eps}-X_{t}}{\eps}> dt.
$$
Whenever they exist we will write
$$
I^\circ (\varphi) = \int_0^T <\varphi(X_t), d^0 X_t>, \quad \text{and} \quad
I^-(\varphi) = \int_0^T <\varphi(X_t), d^- X_t>.
$$

In~\cite{Alos} the authors show that the symmetric integral exists for
fBm with any $H > 1/4$ (the proof  about the $d=1$ case but it
extends without problem to higher dimensions).
For $H > 1/6$ necessary and sufficient conditions for the existence of
the symmetric integral are given in~\cite{Grad}.

As far as the forward integral is concerned, it is known that, in
dimension $1$, it does not exist at least for some (very simple)
functions $\varphi$. Existence of the forward integral in dimension
one is guaranteed
if and only if $H \ge 1/2$~\cite{RV, RVSem}.
Observe that, when $H < \frac{1}{2}$ the forward integral $\int_0^T B d^-B$ 
does not exist since $B$ is not of finite quadratic variation process.

When $H \ge 1/2$ and for $d > 1$ the forward integral (equal to the
Young integral in the case $H > 1/2$), is equal to the symmetric integral minus
the covariation $[\varphi(X),X]/2$. This exists if $X$ has all its \emph{mutual
covariations} and the $i$-th component of that bracket gives
$$[f(X),X^i]_t
=   \int_0^t \sum_{j=1}^d \partial_{j} f (X_s) d[X^i, X^j]_s,
    $$
see again~\cite{RV}.

\medskip

Before entering into details  concerning stochastic integration,
we state an important preliminary result.

\begin{proposition} \label{KAPPA-ALPHA}
If $\alpha > \max(0, d/2-1/(2H))$
 then
$ \expect \left (\int_0^T \int_0^T
 K_\alpha(X_t - X_s) ds dt  \right) < +\infty $.

Therefore Condition~\ref{cond:B} holds.
\end{proposition}
\begin{proof}
Since the corresponding random variable is non-negative, previous expectation equals
\begin{equation} \label{ETech}
  \int_0^T \int_0^T \expect(  K_\alpha(X_t - X_s) )ds dt  =
2  \int_0^T \int_0^t \expect(  K_\alpha(X_t - X_s) )ds dt.
\end{equation}
Next we use the bound on $K_\alpha$ given in Lemma~\ref{lemma su K}:
\begin{itemize}
\item Suppose first that
 $ d/2-1/(2H) <  \alpha < \frac{d}{2}$.
Then the  previous expression
is bounded by 
$
 \const  \int_0^T \int_0^t \expect(   | X_t - X_s |^{2\alpha -d}  )ds dt  
$.
For any $\gamma > -d$, the scaling property  of 
fractional Brownian motion gives 
 $\expect[| X_t - X_s |^\gamma] = \expect [| N |^{\gamma }]
  (t-s)^{\gamma H}$, where $N$ is a standard
$d$-dimensional Gaussian random variable. Then 
the right member of  (\ref{ETech}) equals
\begin{equation} \label{eq:tech} 
 \const\ \expect[\vert N \vert^{2 \alpha   -d} \int_0^T \int_0^t  
  \vert t-s \vert^{2 \alpha H - dH  }  ds dt
\end{equation}
which is  finite if $2 \alpha H  - dH  > -1$ and $\alpha > 0$. Therefore
when $\alpha > \max(d/2-1/(2H),0)$.
\item Suppose now that $\frac{d}{2} < \alpha$.
In this case    $K_\alpha \le \const  $ 
$(\ref{ETech})$ is trivially bounded.
\item Finally, in the case $\alpha = d/2$ we have $K_\alpha(x) \le c
\log |x|$ and  using again the scaling, it  is easy
 to prove the boundedness of  (\ref{ETech}).
\end{itemize}
\end{proof}
\bigskip

The Proposition above will allow us to verify the condition required by 
Lemma \ref{ltech5} of previous section.

\subsection{Symmetric integral}

Let $D^0_\eps X_t = (X_{t+\eps}-X_{t-\eps})/2\eps$ and
let $f:\mathbb{R}^{d} - \{0\}\rightarrow
\mathbb{R}_{+}$ such that
\begin{equation}
  \label{eq:f-class}
\int_0^T \int_0^T f(X_t-X_s) dt ds < \infty\qquad \prob-a.s.  
\end{equation}
 and denote
\[
Z_{\eps}(f):=\int_{0}^{T}\int_{0}^{T}f\left(
X_{t}-X_{s}\right)  
\langle D^0_\eps X_t, D^0_\eps X_s \rangle
dsdt.
\]

\bigskip
In this section we will study the r.v. $Z_{\alpha,\eps} =
Z_\eps(K_\alpha)$ in order to obtain necessary and sufficient
conditions for the regularity of the symmetric current $I(\varphi)$
based on fBm. 

The following lemma is proved in Appendix~A.
\begin{lemma}  
\label{lemma:bounds} 
\begin{enumerate}
\item We have the following estimates:
 \begin{equation}
    \label{eq:bound1}
     \vert\cov(D^0_\eps X^i_t, D^0_\eps X^i_s)| 
\le \const |t-s|^{2H-2};
  \end{equation}
and
\begin{equation}
  \label{eq:bound2}
|\cov(D^0_\eps X^i_t,  X^i_{t}-X^i_{s})|   = |\cov(D^0_\eps X^i_s,  X^i_{t}-X^i_{s})|
\le \const |t-s|^{2H-1}  
\end{equation}
uniformly in $\eps>0$.
\item If $s \neq t$, one has
$$
\lim_{\eps \to 0} |\cov(D^0_\eps X^1_t,X^1_t- X^1_s)| =  2 H  |t-s|^{2H-1}.
$$
\end{enumerate}
\end{lemma}

The main results of this section are contained in the next two theorems. 
Let $\alpha_H = d/2-1+1/(2H)$.

\begin{theorem}
 \label{th:reg-symm}
For any $H>1/4$ and any $\alpha >\alpha_H$ we have
$\sup_{\eps} \expect Z_{\alpha,\eps} < \infty$.
\end{theorem}
\begin{remark}  \label{re:reg-symm}
We observe that $d/2-1/2H < \alpha_H$ so that
for the range value for $\alpha$ in Theorem \ref{th:reg-symm},
Condition \ref{cond:B} is verified.
\end{remark}
\begin{theorem}
\label{th:opt-symm}
For any $H > 1/4$ and any $\alpha <\alpha_H$ 
  we have
$\liminf_{\eps} E(Z_{\alpha,\eps}) = +\infty$.
\end{theorem}
\begin{remark}\label{Ropt}
The statement of  Theorem \ref{th:opt-symm} and Lemma \ref{linter}
below allow us to formulate the following conjecture that we have not been
able to prove: for $H>1/4$ and $\alpha < \alpha_H$ we should have
$$ \liminf_{\eps \rightarrow 0} Z_{\alpha,\eps} = + \infty, \quad \rm{a.s.} $$

For $\{I(\varphi)\}_{\varphi \in V}$ to have a pathwise redefinition
 on some $A \in \mathcal{F}$,$\prob(A) > 0$ a necessary condition is
that $\limsup_{\eps \to 0} Z_{\alpha,\eps} < \infty$ on $A$ (Th.~\ref{th:neg-res-l2}).

If  this conjecture  were true we could establish 
that 
a.s. $I(\varphi)$ does not belong to $H^{-\alpha}_{2}(\RR^d)$ when
$\alpha < \alpha_H$.
\end{remark}

Before proving the Theorems we deduce  Sobolev regularity of fBm
with any Hurst parameter between $(1/4,1)$.

\begin{corollary}[Regularity of symmetric currents]
\label{cor:current-l2-reg}  
The symmetric integral of a fractional Brownian motion with Hurst
parameter $H > 1/4$ admits a pathwise redefinition  on the space
$H^{-\alpha}_{2}(\RR^d)$ for any $\alpha > \alpha_H$.
\end{corollary}
\begin{proof}
By Theorem~\ref{th:reg-symm} we know that $\expect[ Z_{\alpha,\eps}]$ is
uniformly bounded in $\eps$ when $\alpha > \alpha_H$.
Since the regularized integrals $I_\eps(\varphi)$ converges in probability
 as $\eps \to 0$ for any $H > 1/4$, Condition~\ref{cond:A} holds. So
we can   apply
Theorem~\ref{th:l2result} to obtain a pathwise current with values
in $H^{-\alpha}_{2}(\RR^d)$ for any $\alpha > \alpha_H$.
\end{proof}

\bigskip
In the following proof we will use a basic result about Gaussian random variables recalled here:

\begin{lemma}[``Wick's theorem'']
\label{th:wick}
Let $Z = (Z_{\ell})_{1 \le \ell \le N}$ be 
a mean-zero Gaussian random vector and $f \in C^{1}(\R^{N} ;\R)$, then we have
\begin{equation} \label{Ewick}
\expect [Z_{\ell} f(Z)] = \sum_{j=1}^{N} \cov(Z_{\ell},Z_{j})\expect 
[\nabla_{j} f(Z)].
\end{equation}
\end{lemma}
 \begin{proof}
The conclusion follows easily taking 
$f(z_1,...,z_N) = \exp(i \sum_{j=1}^N t_j z_j)$ 
for any  $t=(t_1,...,t_N) \in \R^N $.
In fact,  in that case  $\expect(f(Z)) $ is provided by the
characteristic function. Therefore one has
$$ \expect(\exp(it \cdot Z)) =  \exp \left(- \frac{t \Gamma t'}{2} \right),$$
where $t'$ stands for the transposition and $\Gamma$ is the covariance
matrix of $Z$.
Differentiating the previous expression with respect to $t_\ell$
provides the result (\ref{Ewick}) for the particular case of $f$.
The general result follows by usual density arguments.
\end{proof}

\bigskip
\begin{proof}[Proof of Th.~\ref{th:reg-symm}]
Using Lemma~\ref{th:wick}, independence and equal distribution of different coordinates we have
\begin{equation}
\label{eq:wick1}
  \begin{split}
\expect Z_{\alpha,\eps} & =
\expect     
\int_{0}^{T}\int_{0}^{T}K_\alpha\left(
X_{t}-X_{s}\right) \sum_{i}
\cov(D^0_\eps X^i_t, D^0_\eps X^i_s)
\,dt ds
\\ & \qquad +
\expect     
\int_{0}^{T}\int_{0}^{T} \Delta K_\alpha\left(
X_{t}-X_{s}\right)  
\cov(D^0_\eps X^1_t,  X^1_{t}-X^1_{s})
\cov(D^0_\eps X^1_s,  X^1_{t}-X^1_{s})
\,dt ds
\\  & =
\expect     
\int_{0}^{T}\int_{0}^{T}K_\alpha\left(
X_{t}-X_{s}\right) \sum_{i} 
\cov(D^0_\eps X^i_t, D^0_\eps X^i_s)
\,dt ds
\\ & \qquad + 
\expect     
\int_{0}^{T}\int_{0}^{T} (-\Delta) K_\alpha\left(
X_{t}-X_{s}\right)  
|\cov(D^0_\eps X^1_t,  X^1_{t}-X^1_{s})|^2
\,dt ds
  \end{split}
\end{equation}
where we used the fact that
$$
\cov(D^0_\eps X^1_t,  X^1_{t}-X^1_{s}) = -
\cov(D^0_\eps X^1_s,  X^1_{t}-X^1_{s}) = \frac{1}{2\eps}\left(|t-s+\eps|^{2H}-|t-s-\eps|^{2H}\right)
$$
which can be verified by a straightforward computation.

 Consider first the case $H > 1/2$. Assume $\alpha < d/2$ and note that when $H>1/2$ we have $\alpha_H < d/2$.
By lemma~\ref{lemma:bounds} we have
\begin{equation*}
  \begin{split}
\expect Z_{\alpha,\eps} & \le \const\,
\expect     
\int_{0}^{T}\int_{0}^{T}K_\alpha\left(
X_{t}-X_{s}\right)  |t-s|^{2H-2}
\,dt ds
\\ & \qquad 
+ \const\,
\expect     
\int_{0}^{T}\int_{0}^{T} |(-\Delta) K_\alpha\left(
X_{t}-X_{s}\right)|  |t-s|^{4H-2}
\,dt ds
  \end{split}
\end{equation*}
and then, using lemma \ref{lemma su K} 1) and 4),
we get
\begin{equation}
\label{eq:ubound}
  \begin{split}
\expect Z_{\alpha,\eps} & \le \const\,
\expect     
\int_{0}^{T}\int_{0}^{T} |X_t-X_s|^{2\alpha-d}  |t-s|^{2H-2}
\,dt ds
\\ & \qquad 
+ \const\,
\expect     
\int_{0}^{T}\int_{0}^{T} |X_t-X_s|^{2\alpha-d-2}  |t-s|^{4H-2}
\,dt ds
\\ & = \const\,
\int_{0}^{T}\int_{0}^{T}  |t-s|^{(2\alpha+2-d)H-2}
\,dt ds
\\ & \qquad 
+ \const\,
\int_{0}^{T}\int_{0}^{T} |t-s|^{(2\alpha-d+2)H-2}
\,dt ds
  \end{split}
\end{equation}
which are uniformly bounded in $\eps$ if $(2\alpha+2-d)H > 1$, i.e. when $\alpha > \alpha_H$ as required. We have established the uniform bound when $\alpha \in (\alpha_H, d/2)$ but now recall that $E Z_{\alpha,\eps}$ is a decreasing function of $\alpha$ so that this bound extends to all $\alpha > \alpha_H$.

Let us now consider the case $H \le 1/2$ (so that now $\alpha_H >d/2$) and assume $\alpha > d/2$.
 Rewrite $Z_{\alpha,\eps}$
 as
 \begin{equation}
   \label{eq:Zreg}
Z_{\alpha,\eps} = Z_{\eps}(h_\alpha) +
Z_{\eps}(K_\alpha(0))
 \end{equation}
where $h_\alpha(x) = K_\alpha(x) - K_\alpha(0)$. Note that
$0 < K_\alpha(0) < \infty$ when $\alpha > d/2$.
Moreover
$$
Z_{\eps}(K_\alpha(0)) =  \int_0^T \int_0^T 
K_\alpha(0) \langle D^0_\eps X_s,D^0_\eps X_t \rangle =  
K_\alpha(0)  \left|\int_0^T  D^0_\eps X_s ds\right|^2
$$
and
$$
\int_0^T  D^0_\eps X_t dt = \int_0^T \frac{ X_{t+\eps}- X_{t-\eps}}{\eps}
dt = \eps^{-1} \int^{T+\eps}_{T} X_t dt 
$$
so that by the continuity of the process $X$ we have the limit
$$
\lim_{\eps \to 0} Z_{\eps}(K_\alpha(0)) = K_\alpha(0) |X_T|^2
$$
exists almost surely and $Z_{\eps}(K_\alpha(0))$ is uniformly in
$L^1$ (actually in all $L^p$).

So it remains to consider $Z_{\eps}(h_\alpha)$. For $h_\alpha$ we have
the estimate~(\ref{eq:asympK2}), provided $0 < 2\alpha-d \le 2$, so that we obtain an upper bound
similar to eq.~(\ref{eq:ubound}) and the same condition on $\alpha$
follows. Note that $\alpha$ must satisfy
$$
d/2 <  \alpha_H < \alpha < \frac{d}{2}+1
$$
so that we must require $H > 1/4$. Again by monotonicity of $\alpha \mapsto Z_{\alpha,\eps}$ we have uniform boundedness for any $\alpha > \alpha_H$.

\end{proof}

\begin{proof}[Proof of Th.~\ref{th:opt-symm}]
We will perform a decomposition of $Z_{\alpha,\eps}$ as follows. 
Write
$$
Z_{\alpha,\eps} = A_{\alpha,\eps} + B_{\alpha,\eps} + Q_{\alpha,\eps}
$$
where
$$
A_{\alpha,\eps} = 
\int_{0}^{T}\int_{0}^{T}K_\alpha\left(
X_{t}-X_{s}\right)  \sum_{i}
\cov(D^0_\eps X^i_t, D^0_\eps X^i_s)
\,dt ds
$$
$$
B_{\alpha,\eps} = 
\int_{0}^{T}\int_{0}^{T} (-\Delta) K_\alpha\left(
X_{t}-X_{s}\right)  
|\cov(D^0_\eps X^1_t,  X^1_{t}-X^1_{s})|^2
\,dt ds
$$
and $Q_{\alpha,\eps}$ is the remainder. Note that, by comparing this
decomposition with eq.~(\ref{eq:wick1}), we have
$$
\expect Z_{\alpha,\eps} = \expect A_{\alpha,\eps} + \expect
B_{\alpha,\eps} 
$$
so that $\expect Q_{\alpha,\beta} = 0$. This is a kind of Wick
product decomposition, but not quite, since the terms $A,B$ are not
constants, but still random variables.

A useful remark is that $A_{\alpha,\eps} \ge 0$ since we can write
$$
A_{\alpha,\eps} = \hat \expect
\int_{0}^{T}\int_{0}^{T}K_\alpha\left(
X_{t}-X_{s}\right)  \langle
D^0_\eps \hat X_t, D^0_\eps \hat X_s\rangle
\,dt ds
$$
where we introduced an auxiliary independent $d$-dimensional fBm $\hat X$ with the same distribution 
of $X$ and where $\hat \expect$ denotes expectation with respect to this auxiliary fBm.
So we have the formula
$
A_{\alpha,\eps} = \hat \expect \|\hat \eta_{\eps}\|^2
$
where
$
\hat \eta_{\eps}(x) = \int_0^T K_{\alpha/2}(x-X_t) D_\eps^0 \hat X_t dt
$
which shows that $A_{\alpha,\eps}> 0$.

Next using the equality $-\Delta K_\alpha(x) = K_{\alpha-1}(x)-K_\alpha(x)$ we rewrite $B_{\alpha,\eps}$ as $B^{(1)}_{\alpha,\eps} - B^{(2)}_{\alpha,\eps}$ where
\begin{equation*}
B^{(1)}_{\alpha,\eps} = 
\int_{0}^{T}\int_{0}^{T} K_{\alpha-1}\left(
X_{t}-X_{s}\right)  
|\cov(D^0_\eps X^1_t,  X^1_{t}-X^1_{s})|^2
\,dt ds
\end{equation*}
and
\begin{equation*}
B^{(2)}_{\alpha,\eps} = 
\int_{0}^{T}\int_{0}^{T} K_\alpha(X_t-X_s)   
|\cov(D^0_\eps X^1_t,  X^1_{t}-X^1_{s})|^2
\,dt ds.
\end{equation*}

Let us show first that $E| B^{(2)}_\eps|$ is uniformly bounded in $\eps$ when $\alpha > \alpha_H-1$ and $H> 1/4$. Indeed when $\alpha \le d/2$ by computations similar to those of Th.~\ref{th:reg-symm} we have that 
$E |B^{(2)}_\eps|$ is uniformly bounded if $\alpha > \alpha_H-1$ and
$\alpha > 0$. On the other hand, when $\alpha > d/2$, the kernel
$K_\alpha$ is bounded, so Lemma \ref{lemma:bounds} 1) allows  to write 
\begin{equation*}
E |B^{(2)}_{\alpha,\eps}| \le \const
\int_{0}^{T}\int_{0}^{T} 
|\cov(D^0_\eps X^1_t,  X^1_{t}-X^1_{s})|^2
\,dt ds \le \const \int_{0}^{T}\int_{0}^{T} 
|t-s|^{4H-2}
\,dt ds \le C
\end{equation*}
uniformly in $\eps$ provided $H > 1/4$.
 
Moreover below we will show the following.

\begin{lemma} \label{linter}
If $H > 1/4$ and $\alpha <\alpha_H$ then $\liminf_{\eps \rightarrow 0} B^{(1)}_{\alpha,\eps}  = + \infty $ a.s.
\end{lemma}

Then if we admit the result of previous lemma we can conclude with the use of Fatou lemma.
In fact, for any $\alpha \in (\alpha_H-1,\alpha_H)$ and for some positive constant $c$ we have 
\begin{equation*}
  \begin{split}
\liminf_{\eps \rightarrow 0} \expect (Z_{\alpha,\eps}) & \ge 
 \liminf_{\eps \rightarrow 0} \expect (B_{\alpha,\eps}) \ge \expect (\liminf_{\eps \rightarrow 0} B^{(1)}_{\alpha,\eps}) - \sup_\eps \expect ( B^{(2)}_{\alpha,\eps})  
\\ & \ge
\expect (\liminf_{\eps \rightarrow 0} B^{(1)}_{\alpha,\eps}) - c = + \infty.    
  \end{split}
\end{equation*}
Moreover observe that this is enough since, from Lemma~\ref{lemma:ordering}, we have that  $Z_{\alpha,\eps}$ is a decreasing function of $\alpha$ so the result will hold for any $\alpha < \alpha_H$.
\end{proof}
\bigskip

\begin{proof}[Proof of Lemma~\ref{linter}]

 By Fatou lemma we have
\begin{equation*}
  \begin{split}
\liminf_{\eps \to 0} B^{(1)}_{\alpha,\eps} & \ge 2    
\int_{0}^{T}\int_{0}^{t} K_{\alpha-1}\left(
X_{t}-X_{s}\right) 
\liminf_{\eps \to 0} |\cov(D^0_\eps X^1_t,X^1_t- X^1_s)|^2
\,dt ds
  \end{split}
\end{equation*}
But, when $t \neq s$,
$$
\liminf_{\eps \to 0} |\cov(D^0_\eps X^1_t,X^1_t- X^1_s)|^2 = 4 H^2 |t-s|^{4H-2}
$$.
Now assume \begin{equation}
  \label{eq:cond1-alpha}
\alpha < d/2+1.  
\end{equation}
then by Lemma~\ref{lemma su K} there exist a small constant $r>0$  such that $K_{\alpha-1}(x) \ge C |x|^{2\alpha-d-2}1_{B(0,r)}(x)$. This allows us to bound from below as follows
\begin{equation*}
  \begin{split}
\liminf_{\eps \to 0} B^{(1)}_{\alpha,\eps} & \ge \const\,   
\int_{0}^{T}\int_{0}^{t}  K_{\alpha-1}\left(
X_{t}-X_{s}\right)  |t-s|^{4H-2}
\,dt ds
\\  & \ge \const\,   
\int_{0}^{T}\int_{0}^{t} |X_t-X_s|^{2\alpha-d-2}
1_{B(0,r)}(X_t-X_s)
  |t-s|^{4H-2}
\,dt ds.
  \end{split}
\end{equation*}
Since the paths of fBm are H\"older continuous with parameter strictly smaller than $H$,
 for any $\gamma < H$ there exists a random constant
$C_{X,\gamma}$ such that
$$
|X_t-X_s| \le C_{X,\gamma} |t-s|^\gamma, \qquad t,s \in [0,T]
$$
By choosing a random time $S>0$ small enough such that $\sup_{t,s \in [0,S]} |X_t-X_s| < r$
we have
\begin{equation}
\label{eq:div-B}
  \begin{split}
\liminf_{\eps \to 0}B^{(1)}_{\alpha,\eps}  & \ge \const\,
\int_{0}^{S}\int_{0}^{t} |X_t-X_s|^{2\alpha-d-2}
  |t-s|^{4H-2}
\,dt ds
\\   & \ge \const\, C_{X,\gamma}^{d-2\alpha+2}   
\int_{0}^{S}\int_{0}^{t} 
  |t-s|^{2H(\alpha-\alpha_H)-1+\delta}
\,dt ds 
  \end{split}
\end{equation}
where $\delta = (d+2-2\alpha)(H-\gamma)$ is a arbitrarily small positive
constant since $\gamma < H$ can be chosen arbitrarily near to $H$ and
$$
d+2-2\alpha > 4-\frac{1}{H} > 0
$$
when $H>1/4$.
Then when $\alpha < \alpha_H$
we can choose $\delta$ small enough to make the double integral
in eq.~(\ref{eq:div-B}) diverge. 
Summing up we must have 
$
\alpha < \min (\alpha_H, d/2+1) \quad \text{and} \quad H > 1/4.
$ and $H> 1/4$. But when $H>1/4$ we have $\alpha_H < d/2-1$ so  that sufficient conditions are $\alpha <\alpha_H$ and $H>1/4$.
This observation concludes the proof.
\end{proof}

\subsection{The Forward integral}

Let $D^-_\eps X_t = (X_{t+\eps}-X_{t})/\eps$,
take $f:\mathbb{R}^{d}-\{0\}\rightarrow
\mathbb{R}_{+}$ satisfying~(\ref{eq:f-class}) and denote
\[
\wZ_{\eps}(f)=\int_{0}^{T}\int_{0}^{T}f\left(
X_{t}-X_{s}\right)  
\langle D^-_\eps X_t, D^-_\eps X_s \rangle
dsdt.
\]

We can state similar theorems as in previous subsection.
\begin{theorem}
 \label{th:reg-forw}
For any $H \ge 1/2$ and any $\alpha >\alpha_H$ we have
$\sup_{\eps} \expect Z_{\alpha,\eps} < \infty$.
\end{theorem}

\begin{theorem}
\label{th:opt-forw}
For $H \ge 1/2$ and  $\alpha < \alpha_H$ 
we have $\liminf_{\eps} E(Z_{\alpha,\eps}) = +\infty$.
\end{theorem}

The proofs follow the same line as the corresponding theorems about 
symmetric integrals. For this the following lemma will be crucial.

\begin{lemma}  
\label{lemma:bounds-forw}
\begin{enumerate}
\item We have the following estimates:
  \begin{equation}
    \label{eq:bound1f}
     \vert\cov(D^-_\eps X^i_t, D^-_\eps X^i_s)| 
\le \const |t-s|^{2H-2}
  \end{equation}
and
\begin{equation}
  \label{eq:bound2f}
|\cov(D^-_\eps X^i_t,  X^i_{t}-X^i_{s})|   = |\cov(D^-_\eps X^i_s,  X^i_{t}-X^i_{s})|
\le \const |t-s|^{2H-1}  
\end{equation}
\item If $s \neq t$, one has
$$
\lim_{\eps \to 0} |\cov(D^-_\eps X^1_t,X^1_t- X^1_s)| = 2 H |t-s|^{2H-1}.
$$
\end{enumerate}
\end{lemma}

\begin{proof}[Proof of Th.~\ref{th:opt-forw}]
Again the proof is  similar to the one of  Th.~\ref{th:opt-symm}
where $D^-$ is replaced by $D^0$, Lemma \ref{lemma:bounds-forw} is used
instead of lemma \ref{lemma:bounds}.
In particular, according to Lemma \ref{lemma:bounds-forw} 
 if $s \neq t$ one has
$$
\liminf_{\eps \to 0} |\cov(D^-_\eps X^1_t,X^1_t- X^1_s)|^2 = 4 H^2 |t-s|^{4H-2}
$$
when $H \ge \frac{1}{2}$.
\end{proof}

In particular we have the following Corollary.
\begin{corollary}[Regularity of forward currents]
\label{cor:current-l2-regforw}  
The forward integral of a fractional Brownian motion with Hurst
parameter $H \ge 1/2$ admits a pathwise redefinition  on the space
$H^{-\alpha}_{2}(\RR^d)$ for any $\alpha > \alpha_H$.
\end{corollary}
\begin{proof} 
By Theorem~\ref{th:reg-forw} we know that $\expect Z_{\alpha,\eps}$ is
uniformly bounded in $\eps$ when $\alpha > \alpha_H$.
Since the regularized integrals $I^-_\eps(\varphi)$ converges in probability as
 $\eps \to 0$ for any $H \ge 1/2$, Condition~\ref{cond:A} holds and  we can  apply
Theorem~\ref{th:l2result} to obtain a pathwise current with values
in $H^{-\alpha}_{2}(\RR^d)$ for any $\alpha > \alpha_H$.
\end{proof}

\section{Brownian regularity in $H_{p}^{\alpha}$, $p\neq2$}

In this section we restrict ourselves to the case when $X$ is a $d$%
-dimensional classical Brownian motion, that we denote by $W$. The key
ingredient is the following lemma.

\begin{lemma}
\label{lemma stima base su W}If the dimension $d\geq2$ and the real numbers
$\alpha>1$ and $p^{\prime}>1$ satisfy
\[
\left(  d-\alpha+1\right)  p^{\prime}<d
\]
then
\[
\int_{\mathbb{R}^{d}}E\left[  \left(  \int_{0}^{T}\frac{\exp\left(
-\varepsilon\left|  x-W_{t}\right|  \right)  }{\left|  x-W_{t}\right|
^{2d-2\alpha}}dt\right)  ^{p^{\prime}/2}\right]  dx<\infty
\]
for every $\varepsilon>0$.
\end{lemma}

We shall prove below this lemma. Let us first describe its consequences.

\bigskip From the bounds on $K_{\alpha/2}$, see \ref{lemma su K},  we have
\[
\int_{\mathbb{R}^{d}}E\left[  \left(  \int_{0}^{T}K_{\alpha/2}^{2}\left(
x-W_{t}\right)  dt\right)  ^{p^{\prime}/2}\right]  dx<\infty
\]
and thus for a.e. $x\in\mathbb{R}^{d}$ we have
$
P(  \int_{0}^{T}K_{\alpha/2}^{2}(  x-W_{t})  dt<\infty)
=1
$
which implies that the It\^{o} integral
\[
\eta\left(  x\right)  :=\int_{0}^{T}K_{\alpha/2}\left(  x-W_{t}\right)  dW_{t}%
\]
is well defined, for a.e. $x\in\mathbb{R}^{d}$, as a limit in probability of
\[
\int_{0}^{T}K_{\alpha/2}\left(  x-W_{t}\right)  \frac{W_{t+\varepsilon}-W_{t}%
}{\varepsilon} dt,
\]
see for instance \cite{RVSem}.
These approximation integrals are measurable in the pair $\left(
x,\omega\right)  $, hence they are measurable in $x$ as a mapping with values
in the space of random variables with the metric of convergence in
probability, and this way one can see that the limit object $\eta\left(
x\right)  $ is measurable in the pair $\left(  x,\omega\right)  $. From
Burkhoder-Davies-Gundy (BDG) inequality we have
$
\int_{\mathbb{R}^{d}}E\left[  \left|  \eta\left(  x\right)  \right|
^{p^{\prime}}\right]  dx<\infty
$
and thus $\eta\in L^{p^{\prime}}\left(  \Omega\times\mathbb{R}^{d}\right)  $
and
\[
P\left(  \omega\in\Omega:x\mapsto\eta\left(  x,\omega\right)  \in
L^{p^{\prime}}\left(  \mathbb{R}^{d}\right)  \right)  =1.
\]

To minimize the subtleties related to a direct use of $\eta\left(  x\right)
$, we introduce a regularization. Let $T\left(  t\right)  $ be the semigroup
on $L^{2}\left(  \mathbb{R}^{d}\right)  $ generated by $(\Delta-1)$ and we
set
\[
K_{\alpha/2}^{\left(  \delta\right)  }\left(  x\right)  :=\left(  T\left(
\delta\right)  K_{\alpha/2}\right)  \left(  x\right)  =\left(  4\pi
\delta\right)  ^{-d/2}\int_{\mathbb{R}^{d}}e^{-\frac{\left|  x-y\right|  ^{2}%
}{4\delta}-\delta}K_{\alpha/2}\left(  y\right)  dy.
\]
We have $K_{\alpha/2}^{\left(  \delta\right)  }\in S\left(  \mathbb{R}%
^{d}\right)  $. Set
\[
\eta^{\left(  \delta\right)  }\left(  x\right)  :=\int_{0}^{T}K_{\alpha
/2}^{\left(  \delta\right)  }\left(  x-W_{t}\right)  dW_{t}%
\]
which is obviously well defined for every $x$ and has a measurable version in
the pair $\left(  x,\omega\right)  $. It is not difficult to justify that
$\eta^{\left(  \delta\right)  }$ is square integrable in $\left(
x,\omega\right)  $ and that $\eta^{\left(  \delta\right)  }=T\left(
\delta/2\right)  \eta^{\left(  \delta/2\right)  }$ hence $\eta^{\left(
\delta\right)  }\in S\left(  \mathbb{R}^{d}\right)  $ with probability one.

We have the following regularized version of (\ref{decoupling}).

\begin{lemma}
For $\delta>0$,
\begin{equation}
\int_{0}^{T}\left\langle \left(  T\left(  \delta\right)  \varphi\right)
\left(  W_{t}\right)  ,dW_{t}\right\rangle =\int_{\mathbb{R}^{d}}\left\langle
(1-\Delta)^{\alpha/2}\varphi\left(  x\right)  ,\eta^{\left(  \delta\right)
}\left(  x\right)  \right\rangle dx.\label{decoupling W delta}%
\end{equation}
\end{lemma}

\begin{proof}
Given a vector field $\phi\in S\left(  \mathbb{R}^{d}\right)  $ and a
continuous exponentially decreasing function $\psi$ on $\mathbb{R}^{d}$, we
have the Fubini type identity
\[
\int_{0}^{T}\left\langle \int_{\mathbb{R}^{d}}\psi\left(  W_{t}-x\right)
\phi(x)dx,dW_{t}\right\rangle =\int_{\mathbb{R}^{d}}\left\langle
\phi(x)dx,\int_{0}^{T}\psi\left(  W_{t}-x\right)  dW_{t}\right\rangle
\]
with probability one. We omit the details of the proof. We have
\[
T\left(  \delta\right)  \varphi=T\left(  \delta\right)  (1-\Delta)^{-\alpha
/2}(1-\Delta)^{\alpha/2}\varphi=\int_{\mathbb{R}^{d}}K_{\alpha/2}^{\left(
\delta\right)  }\left(  \cdot-x\right)  \left[  (1-\Delta)^{\alpha/2}%
\varphi\right]  (x)dx
\]
and thus
\[
\int_{0}^{T}\left\langle \left(  T\left(  \delta\right)  \varphi\right)
\left(  W_{t}\right)  ,dW_{t}\right\rangle =\int_{0}^{T}\left\langle
\int_{\mathbb{R}^{d}}K_{\alpha/2}^{\left(  \delta\right)  }\left(
W_{t}-x\right)  \left[  (1-\Delta)^{\alpha/2}\varphi\right]  (x)dx,dW_{t}%
\right\rangle .
\]
Here we can apply the Fubini rule because $T\left(  \delta\right)  \varphi\in
S\left(  \mathbb{R}^{d}\right)  $ for $\delta>0$ and $(1-\Delta)^{\alpha
/2}\varphi\in S\left(  \mathbb{R}^{d}\right)  $. This implies (\ref{decoupling
W delta}) and completes the proof.
\end{proof}

\begin{lemma}
For $d\geq2$, $\alpha>1$, $p^{\prime}>1$, such that
$
\left(  d-\alpha+1\right)  p^{\prime}<d
$
we have
\begin{equation}
\sup_{\delta>0}\int_{\mathbb{R}^{d}}E\left[  \left|  \eta^{\left(
\delta\right)  }\left(  x\right)  \right|  ^{p^{\prime}}\right]
dx<\infty.\label{stima  base su eta delta}%
\end{equation}
\end{lemma}

\begin{proof}
Let us restrict the argument to the most difficult case $0<\alpha<d$ where
$K_{\alpha/2}$ has a singularity at zero. From BDG inequality we have
\[
\int_{\mathbb{R}^{d}}E\left[  \left|  \eta^{\left(  \delta\right)  }\left(
x\right)  \right|  ^{p^{\prime}}\right]  dx\leq C_{p^{\prime}}\int
_{\mathbb{R}^{d}}E\left[  \left(  \int_{0}^{T}\left|  K_{\alpha/2}^{\left(
\delta\right)  }\left(  x-W_{t}\right)  \right|  ^{2}dt\right)  ^{p^{\prime
}/2}\right]  dx.
\]
From the definition of $K_{\alpha/2}^{\left(  \delta\right)  }$ in terms of
$K_{\alpha/2}$ and estimate (\ref{eq:asympK}) we get the inequality
\[
K_{\alpha/2}^{\left(  \delta\right)  }\left(  x\right)  \leq C_{\alpha
,d}\left(  4\pi\delta\right)  ^{-d/2}\int_{\mathbb{R}^{d}}e^{-\frac{\left|
x-y\right|  ^{2}}{4\delta}}\left|  y\right|  ^{\alpha-d}e^{-\frac{\left|
y\right|  }{8}}dy.
\]
Let us show that this implies%
\begin{equation}
K_{\alpha/2}^{\left(  \delta\right)  }\left(  x\right)  \leq C_{\alpha
,d}\left|  x\right|  ^{\alpha-d}e^{-\frac{\left|  x\right|  }{8}%
}\label{stima su K delta}%
\end{equation}
for a new constant $C_{\alpha,d}$, uniformly in $\delta\in\left(  0,1\right)
$. The proof of this result for $\left|  x\right|  >1$ is rather easy, so let
us only deal with $\left|  x\right|  \leq1$. Write $x=re$ with $\left|
e\right|  =1$ and change variable $y=rz$ in the integral, to get
\begin{align*}
K_{\alpha/2}^{\left(  \delta\right)  }\left(  x\right)  \leq r^{\alpha
-d}C_{\alpha,d}\left(  4\pi\delta\right)  ^{-d/2}r^{d}\int_{\mathbb{R}^{d}%
}e^{-\frac{r^{2}\left|  e-z\right|  ^{2}}{4\delta}}\left|  z\right|
^{\alpha-d}dz\\
=r^{\alpha-d}C_{\alpha,d}\left[  T\left(  \frac{\delta}{r^{2}}\right)  \left|
.\right|  ^{\alpha-d}\right]  \left(  e\right)
\end{align*}
(see the definition of the semigroup $T(t)$). It is now easy to see that
$\left[  T\left(  t\right)  \left|  .\right|  ^{\alpha-d}\right]  \left(
e\right)  $ is bounded above by a constant, uniformly in $t\geq0$. This proves
(\ref{stima su K delta}).

Having this estimate, it is sufficient to apply lemma \ref{lemma stima base su
W}. The proof is complete.
\end{proof}

We can now prove the main result of this section.

\begin{theorem}
The It\^{o} integral $\int_{0}^{T}\left\langle \varphi\left(  W_{t}\right)
,dW_{t}\right\rangle $ has a pathwise redefinition on the space $V=H_{p}%
^{\alpha}\left(  \mathbb{R}^{d}\right)  $ for every dimension $d$ and real
numbers $\alpha>1$ and $p>1$ satisfying
\[
p>\frac{d}{\alpha-1}.
\]
In particular, in any dimension $d$, given $\varepsilon>0$, for every
$p>\frac{d}{\varepsilon}$ the integral $\int_{0}^{T}\left\langle
\varphi\left(  W_{t}\right)  ,dW_{t}\right\rangle $ has a pathwise
redefinition on the space $H_{p}^{1+\varepsilon}\left(  \mathbb{R}^{d}\right)
$.
\end{theorem}

\begin{proof}
\textbf{Step 1}. In the case $d=1$ we have $H_{p}^{\alpha}\left(
\mathbb{R}\right)  \subset C^{1}\left(  \mathbb{R}\right)  $ by Sobolev
embedding theorem (see \cite{Triebel}, section 2.8.1, remark 2). Thus
\[
\int_{0}^{T}\varphi\left(  W_{t}\right)  dW_{t}=\Phi\left(  W_{T}\right)
-\Phi\left(  0\right)  -\frac{1}{2}\int_{0}^{T}\varphi^{\prime}\left(
W_{t}\right)  dt
\]
where $\Phi^{\prime}=\varphi$. This implies the result. We restrict now to the
case $d\geq2$.

\textbf{Step 2}. We pass to the limit in (\ref{decoupling W delta}). Let
us treat the left-hand-side. With easy manipulations we see that
\[
\left(  T\left(  \delta\right)  \varphi\right)  \left(  x\right)  =\left(
2\pi\right)  ^{-d/2}e^{-\delta}\int_{\mathbb{R}^{d}}e^{-\frac{\left|
z\right|  ^{2}}{2}}\varphi\left(  x-z\sqrt{2\delta}\right)  dz
\]
Hence, splitting the integral in a sufficiently large ball and the
complementary, since $\varphi\in S\left(  \mathbb{R}^{d}\right)  $, we see
that $T\left(  \delta\right)  \varphi\rightarrow\varphi$ uniformly over all
$\mathbb{R}^{d}$ as $\delta\rightarrow0$. Thus $\int_{0}^{T}\left\langle
\left(  T\left(  \delta\right)  \varphi\right)  \left(  W_{t}\right)
,dW_{t}\right\rangle $ easily converges to $\int_{0}^{T}\left\langle
\varphi\left(  W_{t}\right)  ,dW_{t}\right\rangle $, in mean square.

Given the value of $p$ in the statement of the theorem, under the assumption
$\alpha>1$ the inequality $p>\frac{d}{\alpha-1}$ is equivalent to $\left(
d-\alpha+1\right)  p^{\prime}<d$, where $1/p + 1/p^\prime = 1$,
so the previous lemma applies. From
(\ref{stima base su eta delta}) there is a sequence $\delta_{n}\rightarrow0$
and an element $\eta^{\left(  0\right)  }\in L^{p^{\prime}}\left(
\Omega\times\mathbb{R}^{d}\right)  $ such that $\eta^{\left(  \delta_n\right)
}\rightharpoonup\eta^{\left(  0\right)  }$ weakly in $L^{p^{\prime}}\left(
\Omega\times\mathbb{R}^{d}\right)  $, when $n \rightarrow +\infty$. From (\ref{decoupling W delta}), for a
given $\varphi\in S\left(  \mathbb{R}^{d}\right)  $, we thus have, in the
limit as $n \rightarrow \infty$,
\[
E\left[  X\int_{0}^{T}\left\langle \varphi\left(  W_{t}\right)  ,dW_{t}%
\right\rangle \right]  =E\left[  X\int_{\mathbb{R}^{d}}\left\langle
(1-\Delta)^{\alpha/2}\varphi\left(  x\right)  ,\eta^{\left(  0\right)
}\left(  x\right)  \right\rangle dx\right]
\]
for every bounded r.v. $X$ and thus
\[
\int_{0}^{T}\left\langle \varphi\left(  W_{t}\right)  ,dW_{t}\right\rangle
=\int_{\mathbb{R}^{d}}\left\langle (1-\Delta)^{\alpha/2}\varphi\left(
x\right)  ,\eta^{\left(  0\right)  }\left(  x\right)  \right\rangle dx
\]
with probability one.

\textbf{Step 3}. Therefore, given $\varphi\in S\left(  \mathbb{R}^{d}\right)
$, with probability one we have
\begin{align*}
\left|  \int_{0}^{T}\left\langle \varphi\left(  W_{t}\right)  ,dW_{t}%
\right\rangle \right|  \leq\left\|  (1-\Delta)^{\alpha/2}\varphi\right\|
_{L^{p}\left(  \mathbb{R}^{d}\right)  }\left\|  \eta^{\left(  0\right)
}\right\|  _{L^{p^{\prime}}\left(  \mathbb{R}^{d}\right)  }\\
\leq C\left\|  \varphi\right\|  _{H_{p}^{\alpha}\left(  \mathbb{R}^{d}\right)
}\left\|  \eta^{\left(  0\right)  }\right\|  _{L^{p^{\prime}}\left(
\mathbb{R}^{d}\right)  }.
\end{align*}
The proof is complete.
\end{proof}

\begin{remark}
\label{remark stopping} The same result is true for the stopped Brownian
motion
\[
W_{t}^{R}=W_{t\wedge\tau_{R}},\quad\tau_{R}=\inf\left\{  t>0:\left|
W_{t}\right|  \geq R\right\}  .
\]
with given $R>0$. The statement is that the It\^{o} integral $\int_{0}%
^{T}\left\langle \varphi\left(  W_{t}^{R}\right)  ,dW_{t}^{R}\right\rangle $
has a pathwise redefinition on the space $H_{p}^{\alpha}\left(  \mathbb{R}%
^{d}\right)  $ under the same conditions on $d,\alpha,p$ as in the theorem.
The proof is the same (even easier, since in the proof of lemma \ref{lemma
stima base su W} we do not have to care of the exponential term).
\end{remark}

\begin{remark}
\label{remark Walphap} The same result is true in the Sobolev-Slobodeckij
spaces $W_{p}^{\alpha}\left(  \mathbb{R}^{d}\right)  $ defined in
\cite{Triebel}, section 2.3. The statement is that the It\^{o} integral
$\int_{0}^{T}\left\langle \varphi\left(  W_{t}\right)  ,dW_{t}\right\rangle $
has a pathwise redefinition on the space $W_{p}^{\alpha}\left(  \mathbb{R}%
^{d}\right)  $ under the same conditions on $d,\alpha,p$ as in the theorem.
Indeed, given a triple $d,\alpha,p$ as in the theorem, let $\alpha^{\prime
}<\alpha$ be such that also the triple $d,\alpha^{\prime},p$ satisfies the
assumption of the theorem. Then $\int_{0}^{T}\left\langle \varphi\left(
W_{t}\right)  ,dW_{t}\right\rangle $ has a pathwise redefinition on
$H_{p}^{\alpha^{\prime}}\left(  \mathbb{R}^{d}\right)  $; by definition of
pathwise redefinition, we see that this implies that $\int_{0}^{T}\left\langle
\varphi\left(  W_{t}\right)  ,dW_{t}\right\rangle $ has a pathwise
redefinition on the space $W_{p}^{\alpha}\left(  \mathbb{R}^{d}\right)  $,
because we have the continuous embedding
\[
W_{p}^{\alpha}\left(  \mathbb{R}^{d}\right)  \subset H_{p}^{\alpha^{\prime}%
}\left(  \mathbb{R}^{d}\right),
\]
see \cite{Triebel}, remark 4 of section 2.3.3. The same result is of course
true for the It\^{o} integral $\int_{0}^{T}\left\langle \varphi\left(
W_{t}^{R}\right)  ,dW_{t}^{R}\right\rangle $.
\end{remark}

We can now elaborate the previous results in the direction of the H\"{o}lder
topology. Given $\varepsilon\in\left(  0,1\right)  $, denote by
$C^{1+\varepsilon}\left(  \mathbb{R}^{d}\right)  $ the space of all
continuously differentiable functions $f$ on $\mathbb{R}^{d}$\ such that
\[
\left\|  f\right\|  _{C^{1+\varepsilon}}=\sup_{x\in\mathbb{R}^{d}}\left(
\left|  f(x)\right|  +\left|  Df(x)\right|  \right)  +\sup_{x\neq y}%
\frac{\left|  Df(x)-Df(y)\right|  }{\left|  x-y\right|  ^{\varepsilon}}%
<\infty,
\]
see \cite{Triebel}, section 2.7. Endowed with the norm $\left\|  .\right\|
_{C^{1+\varepsilon}}$, the space $C^{1+\varepsilon}\left(  \mathbb{R}%
^{d}\right)  $ is a Banach space.

\begin{theorem}
In any dimension $d$, for every $\varepsilon\in\left(  0,1\right)  $ the
It\^{o} integral $\int_{0}^{T}\left\langle \varphi\left(  W_{t}\right)
,dW_{t}\right\rangle $ has a pathwise redefinition on the space
$C^{1+\varepsilon}\left(  \mathbb{R}^{d}\right)  $.
\end{theorem}

\begin{proof}
\textbf{Step 1}. This preliminary step is devoted to a few details used below.
Recall that the classical Sobolev space $W_{p}^{1}\left(  \mathbb{R}%
^{d}\right)  $ is defined as the space of all $f\in L^{p}\left(
\mathbb{R}^{d}\right)  $ having distributional derivative $Df\in L^{p}\left(
\mathbb{R}^{d \times d }  \right)  $. Recall also (see remark 4 of section 2.5.1 of
\cite{Triebel}) that, for every $\varepsilon\in\left(  0,1\right)  $, the
space $W_{p}^{1+\varepsilon}\left(  \mathbb{R}^{d}\right)  $ of remark
\ref{remark Walphap} is characterized as the space of all $f\in W_{p}%
^{1}\left(  \mathbb{R}^{d}\right)  $ such that
\[
\int_{\mathbb{R}^{d}\times\mathbb{R}^{d}}\frac{\left|  Df(x)-Df(y)\right|
^{p}}{\left|  x-y\right|  ^{d+\varepsilon p}}dxdy<\infty
\]
and as a norm on $W_{p}^{1+\varepsilon}\left(  \mathbb{R}^{d}\right)  $ one
can take the following one:%
\[
\left\|  f\right\|  _{W_{p}^{1+\varepsilon}}^{p}=\left\|  f\right\|  _{L_{p}%
}^{p}+\left\|  Df\right\|  _{L_{p}}^{p}+\int_{\mathbb{R}^{d}\times
\mathbb{R}^{d}}\frac{\left|  Df(x)-Df(y)\right|  ^{p}}{\left|  x-y\right|
^{d+\varepsilon p}}dxdy.
\]
Then it is easy to verify that for every $\varepsilon,\varepsilon^{\prime}%
\in\left(  0,1\right)  $ with $\varepsilon>\varepsilon^{\prime}$ the following
assertion is true, where $B\left(  0,R\right)  $ denotes the ball of center
$0$ and radius $R>0$:
\[
f\in C^{1+\varepsilon}\left(  \mathbb{R}^{d}\right)  \text{, }f\text{ with
support in }B\left(  0,R\right)  \Rightarrow f\in W_{p}^{1+\varepsilon
^{\prime}}\left(  \mathbb{R}^{d}\right)
\]
and%
\begin{equation} \label{Esob}
\left\|  f\right\|  _{W_{p}^{1+\varepsilon^{\prime}}}^{p}\leq C\left(
R,\varepsilon,\varepsilon^{\prime},p,d\right)  \left\|  f\right\|
_{C^{1+\varepsilon}}^{p}%
\end{equation}
where $C\left(  R,\varepsilon,\varepsilon^{\prime},p,d\right)  $ is a constant
depending only on $R,\varepsilon,\varepsilon^{\prime},p,d$.

Indeed, we have
\begin{align*}
\int_{\mathbb{R}^{d}\times\mathbb{R}^{d}}\frac{\left|  Df(x)-Df(y)\right|
^{p}}{\left|  x-y\right|  ^{d+\varepsilon^{\prime}p}}dxdy &  =\int_{B\left(
0,R\right)  \times B\left(  0,R\right)  }\frac{\left|  Df(x)-Df(y)\right|
^{p}}{\left|  x-y\right|  ^{d+\varepsilon^{\prime}p}}dxdy\\
&  \leq\int_{\left|  x-y\right|  \leq1,\left|  x\right|  \leq R,\left|
y\right|  \leq R}\frac{\left|  Df(x)-Df(y)\right|  ^{p}}{\left|  x-y\right|
^{d+\varepsilon^{\prime}p}}dxdy \\
&  + \int_{\left|  x-y\right|  >1,\left|  x\right|
\leq R,\left|  y\right|  \leq R}\frac{\left|  Df(x)-Df(y)\right|  ^{p}%
}{\left|  x-y\right|  ^{d+\varepsilon^{\prime}p}}dxdy\\
&  \leq\int_{\left|  x-y\right|  \leq1,\left|  x\right|  \leq R,\left|
y\right|  \leq R}\frac{\left\|  f\right\|  _{C^{1+\varepsilon}}^{p}}{\left|
x-y\right|  ^{d+\left(  \varepsilon^{\prime}-\varepsilon\right)  p}%
}dxdy+C\left(  p,d\right)  \left\|  f\right\|  _{C^{1+\varepsilon}}^{p}R^{d}%
\end{align*}
where $C\left(  p,d\right)  $ is a constant depending only on $p,d$. The claim
(\ref{Esob}) easily follows from this inequality.

\textbf{Step 2}. Let $d$ and $\varepsilon$ be given, as in the claim of the
theorem. Choose $\varepsilon^{\prime}\in\left(  0,\varepsilon\right)  $ and
$p>\frac{d}{\varepsilon^{\prime}}$. Let $W$ be defined on a complete
probability space $\left(  \Omega,\mathcal{F},P\right)  $. Remark \ref{remark
Walphap} states that there exists a random variable $C>0$ such that, for every
$\varphi\in W_{p}^{1+\varepsilon^{\prime}}\left(  \mathbb{R}^{d}\right)  $,
\[
\left|  \int_{0}^{T}\left\langle \varphi\left(  W_{t}\right)  ,dW_{t}%
\right\rangle \right|  \leq C\left\|  \varphi\right\|  _{W_{p}^{1+\varepsilon
^{\prime}}}%
\]
on a full probability set $\Omega_{\varphi}$.

For every $R>0$, let $\theta_{R}:\mathbb{R}^{d}\rightarrow\lbrack0,\infty)$ be
a $C^{\infty}$ function such that $\theta_{R}\left(  x\right)  =1$ for
$\left|  x\right|  \leq R+1$, $\theta_{R}\left(  x\right)  =0$ for $\left|
x\right|  \geq R+2$. Given $\varphi\in C^{1+\varepsilon}\left(  \mathbb{R}%
^{d}\right)  $, we have $\varphi\cdot\theta_{R}\in C^{1+\varepsilon}\left(
\mathbb{R}^{d}\right)  $ and thus $\varphi\cdot\theta_{R}\in W_{p}%
^{1+\varepsilon^{\prime}}\left(  \mathbb{R}^{d}\right)  $. Therefore%
\[
\left|  \int_{0}^{T}\left\langle \left(  \varphi\cdot\theta_{R}\right)
\left(  W_{t}\right)  ,dW_{t}\right\rangle \right|  \leq C\left\|
\varphi\cdot\theta_{R}\right\|  _{W_{p}^{1+\varepsilon^{\prime}}}%
\]
on a full probability set $\Omega_{\varphi\cdot\theta_{R}}$.

From step 1, there exists a random variable $C_{R}>0$, independent of
$\varphi$, such that
\[
\left|  \int_{0}^{T}\left\langle \left(  \varphi\cdot\theta_{R}\right)
\left(  W_{t}\right)  ,dW_{t}\right\rangle \right|  \leq C_{R}\left\|
\varphi\right\|  _{C^{1+\varepsilon}}\text{ on }\Omega_{\varphi\cdot\theta
_{R}},
\]
where we have also used the fact that $\left\|  \varphi\cdot\theta_{R}\right\|
_{C^{1+\varepsilon}}\leq C_{\theta}\left\|  \varphi\right\|
_{C^{1+\varepsilon}}$ for some constant $C_{\theta}>0$ depending on the
function $\theta$ (and thus on $R$ again). Redefine, if necessary, $C_{R}$ in
such a way that $R\mapsto C_{R}$ is non decreasing, with probability one.

Let $A_{R}$ be the set%
\[
A_{R}=\left\{  \tau_{R}>T\right\},
\]
where $\tau_{R}$ is defined in remark \ref{remark stopping}. The sets $A_{R}$
increase with $R$. Given the family of events $A_{R}$ and random variables
$C_{R}$, we can define a new random variable $C^{\prime}>0$ such that
$C_{R}\leq C^{\prime}$ on $A_{R}$ (it is sufficient to put $C^{\prime}%
=C_{N+1}$ on $A_{N+1}\diagdown A_{N}$). Thus, given $\varphi\in
C^{1+\varepsilon}\left(  \mathbb{R}^{d}\right)  $, we have%
\[
\left|  \int_{0}^{T}\left\langle \left(  \varphi\cdot\theta_{R}\right)
\left(  W_{t}\right)  ,dW_{t}\right\rangle \right|  \leq C^{\prime}\left\|
\varphi\right\|  _{C^{1+\varepsilon}}\text{ on }\Omega_{\varphi\cdot\theta
_{R}}\cap A_{R}.
\]

For every $R>0$ and $\varphi\in C^{1+\varepsilon}\left(  \mathbb{R}%
^{d}\right)  $, there is a $P$-null set $N_{R,\varphi}$ such that
\[
\int_{0}^{T}\left\langle \left(  \varphi\cdot\theta_{R}\right)  \left(
W_{t}\right)  ,dW_{t}\right\rangle =\int_{0}^{T}\left\langle \varphi\left(
W_{t}\right)  ,dW_{t}\right\rangle \text{ on }A_{R}\diagdown N_{R,\varphi}.
\]
Therefore, given $R>0$ and $\varphi\in C^{1+\varepsilon}\left(  \mathbb{R}%
^{d}\right)  $, we have%
\[
\left|  \int_{0}^{T}\left\langle \varphi\left(  W_{t}\right)  ,dW_{t}%
\right\rangle \right|  \leq C^{\prime}\left\|  \varphi\right\|
_{C^{1+\varepsilon}}\text{ on }\Omega_{\varphi\cdot\theta_{R}}\cap
A_{R}\diagdown N_{R,\varphi}.
\]
It follows that%
\[
\left|  \int_{0}^{T}\left\langle \varphi\left(  W_{t}\right)  ,dW_{t}%
\right\rangle \right|  \leq C^{\prime}\left\|  \varphi\right\|
_{C^{1+\varepsilon}}\text{ on }\bigcup_{R>0}\left(  \Omega_{\varphi\cdot
\theta_{R}}\cap A_{R}\diagdown N_{R,\varphi}\right)  .
\]
Since $P\left(  \bigcup_{R>0}A_{R}\right)  =1$ we have $P\left(  \bigcup
_{R>0}\left(  \Omega_{\varphi\cdot\theta_{R}}\cap A_{R}\diagdown N_{R,\varphi
}\right)  \right)  =1$. This means $\int_{0}^{T}\left\langle \varphi\left(
W_{t}\right)  ,dW_{t}\right\rangle $ has a pathwise redefinition on the space
$C^{1+\varepsilon}\left(  \mathbb{R}^{d}\right)  $. The proof is complete.
\end{proof}

\begin{remark}
The strategy of step 2 in the previous proof can be used to deal with function
spaces of Fr\'{e}chet type that are not Banach spaces: by localization of the
stochastic process, one can restrict the attention to compact support test
functions and then prove the existence of a pathwise redefinition in
topologies without decay at infinity. For this reason, even the uniformity in
$x\in\mathbb{R}^{d}$ in the definition of $C^{1+\varepsilon}\left(
\mathbb{R}^{d}\right)  $ is not necessary.
\end{remark}

\begin{remark} \label{Rrough}
In rough path theory (see \cite{Lyons}), for every rough path $\gamma$ of a
certain class which includes a.e. path of Brownian motion, a notion of
integral $\int_{0}^{T}\left\langle \varphi\left(  \gamma_{t}\right)
,d\gamma_{t}\right\rangle $ is defined for every function $\varphi$ with
$\varepsilon$-H\"{o}lder first derivative (for arbitrary $\varepsilon>0$). The
previous theorem is conceptually similar; a closer comparison, however,
requires further investigation.
\end{remark}

Finally, we have to prove lemma \ref{lemma stima base su W}.

\subsection{Proof of lemma \ref{lemma stima base su W}}

If $\max_{\left[  0,T\right]  }\left|  W_{t}\right|  \leq\left|  x\right|  /2$
then, for every $t\in\left[  0,T\right]  $,
\begin{align*}
\frac{1}{2}\left|  x\right|    & \leq\left|  x-W_{t}\right|  \leq\frac{3}%
{2}\left|  x\right|  \\
\exp\left(  -\varepsilon\left|  x-W_{t}\right|  \right)    & \leq\exp\left(
-\varepsilon\left|  x\right|  /2\right)  \\
\frac{1}{\left|  x-W_{t}\right|  ^{2d-2\alpha}}  & \leq\frac{\left(
2/3\right)  ^{2d-2\alpha}}{\left|  x\right|  ^{2d-2\alpha}}\text{ \quad if
}2d-2\alpha\leq0\\
\frac{1}{\left|  x-W_{t}\right|  ^{2d-2\alpha}}  & \leq\frac{2^{2d-2\alpha}%
}{\left|  x\right|  ^{2d-2\alpha}}\text{ \quad if }2d-2\alpha>0
\end{align*}
and thus%
\[
\frac{\exp\left(  -\varepsilon\left|  x-W_{t}\right|  \right)  }{\left|
x-W_{t}\right|  ^{2d-2\alpha}}\leq\exp\left(  -\varepsilon\left|  x\right|
/2\right)  \frac{C_{\alpha,d}}{\left|  x\right|  ^{2d-2\alpha}}%
\]
for a suitable constant $C_{\alpha,d}>0$.  Therefore%
\[
\int_{\mathbb{R}^{d}}E\left[  \left(  \int_{0}^{T}\frac{\exp\left(
-\varepsilon\left|  x-W_{t}\right|  \right)  }{\left|  x-W_{t}\right|
^{2d-2\alpha}}dt\right)  ^{p^{\prime}/2}\right]  dx\leq I_{1}+I_{2}%
\]%
\[
I_{1}:=\int_{\mathbb{R}^{d}}\exp\left(  -\varepsilon p^{\prime}\left|
x\right|  /4\right)  E\left[  \left(  \int_{0}^{T}\frac{C_{\alpha,d}}{\left|
x\right|  ^{2d-2\alpha}}dt\right)  ^{p^{\prime}/2}\right]  dx
\]%
$$
I_{2}    :=\int_{\mathbb{R}^{d}}E\left[  1_{\max_{\left[  0,T\right]
}\left|  W_{t}\right|  >\left|  x\right|  /2}\left(  \int_{0}^{T}\frac
{1}{\left|  x-W_{t}\right|  ^{2d-2\alpha}}dt\right)  ^{p^{\prime}/2}\right]
dx.
$$

Obviously $I_{1}<\infty$, being $\left(  d-\alpha+1\right)  p^{\prime}<d$.
Moreover%
\[
I_{2}\leq\int_{\mathbb{R}^{d}}P\left(  \max_{\left[  0,T\right]  }\left|
W_{t}\right|  >\left|  x\right|  /2\right)  ^{\frac{\delta}{1+\delta}}E\left[
\left(  \int_{0}^{T}\frac{1}{\left|  x-W_{t}\right|  ^{2d-2\alpha}}dt\right)
^{\left(  1+\delta\right)  p^{\prime}/2}\right]  ^{\frac{1}{1+\delta}}dx
\]
for every $\delta>0$. Recall the exponential inequality (see \cite{ry}
Proposition 1.8)
\[
P\left(  \max_{t\in\left[  0,T\right]  }W_{t}\geq\beta\right)  \leq
e^{-\frac{\beta}{2T}}.
\]
It easily implies, by symmetry, that
\[
P\left(  \max_{t\in\left[  0,T\right]  }\left|  W_{t}\right|  \geq\left|
x\right|  /2\right)  \leq2e^{-\frac{\left|  x\right|  }{4T}}%
\]
and thus there exist $C_{\delta},\lambda_{\delta}>0$ (depending also on $T$)
such that
\[
P\left(  \max_{t\in\left[  0,T\right]  }\left|  W_{t}\right|  \geq\left|
x\right|  /2\right)  ^{\frac{\delta}{1+\delta}}\leq C_{\delta}e^{-\lambda
_{\delta}\left|  x\right|  }.
\]
Moreover, by Young inequality,
\[
a^{\frac{1}{1+\delta}}=a^{\frac{1}{1+\delta}}\cdot1\leq a+C_{\delta}%
\]
for some constant $C_{\delta}>0$. Thus, for every $\delta>0$,
\[
I_{2}\leq C_{\delta}^{\prime}+C_{\delta}^{\prime}\int_{\mathbb{R}^{d}%
}e^{-\lambda_{\delta}\left|  x\right|  }E\left[  \left(  \int_{0}^{T}\frac
{1}{\left|  x-W_{t}\right|  ^{2d-2\alpha}}dt\right)  ^{\left(  1+\delta
\right)  p^{\prime}/2}\right]  dx
\]
for some constant $C_{\delta}^{\prime}>0$.

The following lemma is inspired from the proof of Corollary 2.4 of Elworthy,
Li, Yor \cite{ELY} and in fact it was suggested to us by K. D. Elworthy.

\begin{lemma}
For every $d\geq2$, $q>1$, $\theta\in\mathbb{R}$, $x\in\mathbb{R}^{d}$, we
have
\begin{align*}
E\left[  \left(  \int_{0}^{T}\frac{dt}{\left|  x+W_{t}\right|  ^{2\left(
1-\theta\right)  }}\right)  ^{\frac q2}\right]  
\leq c_{q,\theta,T}\left\{  E\left[  \left|  x+W_{T}\right|  ^{\theta
q}\right]  +\left|  x\right|  ^{\theta q}+\int_{0}^{T}E\left[  \frac
{dt}{\left|  x+W_{t}\right|  ^{\left(  2-\theta\right)  q}}\right]  \right\}
.
\end{align*}
\end{lemma}

\begin{proof}
Consider the process $Z_{t}=\left|  x+W_{t}\right|  ^{2}$ (squared Bessel
process of dimension $d$). From It\^{o} formula we have
\[
dZ_{t}=2\left\langle x+W_{t},dW_{t}\right\rangle +dt,\quad Z_{0}=\left|
x\right|  ^{2}.
\]
Introducing an auxiliary one-dimensional Brownian motion $\left(  \beta
_{t}\right)  $ we may also write
\[
dZ_{t}=2\sqrt{Z_{t}}d\beta_{t}+dt.
\]
Since $d\geq2$, the one-point sets are polar sets for a $d$-dimensional
Brownian motion, see Proposition 2.7, p. 191 of \cite{ry}. Therefore
$P\left\{  Z_{t}>0,t\in\left[  0,T\right]  \right\}  =1$ and we can develop
$Z_{t}^{\theta/2}$ using It\^{o} formula for any $\theta\in\mathbb{R}$. We
obtain
\begin{align*}
d\left(  Z_{t}^{\theta/2}\right)  =\frac{\theta}{2}Z_{t}^{\left(
\theta-2\right)  /2}\left(  2\sqrt{Z_{t}}d\beta_{t}+dt\right)  +\frac{1}%
{2}\frac{\theta}{2}\frac{\theta-2}{2}Z_{t}^{\left(  \theta-4\right)  /2}%
4Z_{t}dt\\
=\theta Z_{t}^{\left(  \theta-1\right)  /2}d\beta_{t}+c_{\theta}Z_{t}^{\left(
\theta-2\right)  /2}dt
\end{align*}
where $c_{\theta}=\frac{\theta\left(  \theta-1\right)  }{2}$. Therefore
\[
\int_{0}^{T}\theta Z_{t}^{\left(  \theta-1\right)  /2}d\beta_{t}=Z_{T}%
^{\theta/2}-Z_{0}^{\theta/2}-\int_{0}^{T}c_{\theta}Z_{t}^{\left(
\theta-2\right)  /2}dt
\]
and thus, from BDG inequality, for every $q>1$%
\begin{align*}
E\left[  \left(  \int_{0}^{T}\theta^{2}Z_{t}^{\theta-1}dt\right)
^{q/2}\right]  
& \leq c_{q}E\left[  \left(  \int_{0}^{T}\theta Z_{t}^{\left(  \theta-1\right)
/2}d\beta_{t}\right)  ^{q}\right]  \\
& \leq c_{q,\theta}\left\{  E\left[  Z_{T}^{\theta q/2}\right]  +E\left[
Z_{0}^{\theta q/2}\right]  +E\left[  \left(  \int_{0}^{T}Z_{t}^{\left(
\theta-2\right)  /2}dt\right)  ^{q}\right]  \right\}  .
\end{align*}
This implies, by H\"{o}lder inequality,
\begin{align*}
E\left[  \left(  \int_{0}^{T}\frac{1}{Z_{t}^{1-\theta}}dt\right)
^{q/2}\right]  
\leq c_{q,\theta,T}\left\{  E\left[  Z_{T}^{\theta q/2}\right]  +E\left[
Z_{0}^{\theta q/2}\right]  +E\left[  \int_{0}^{T}\frac{1}{Z_{t}^{\left(
2-\theta\right)  q/2}}dt\right]  \right\}
\end{align*}
and the proof is complete.
\end{proof}

\bigskip
We go on with the proof of lemma \ref{lemma stima base su W}.
Simply by taking $1-\theta=d-\alpha$ and $q=\left(  1+\delta\right)
p^{\prime}$ we have:%
\begin{align*}
 &  E\left[  \left(  \int_{0}^{T}\frac{1}{\left|  x-W_{t}\right|  ^{2d-2\alpha
}}dt\right)  ^{\left(  1+\delta\right)  p^{\prime}/2}\right] 
\\ &\qquad \qquad \leq C\left\{  E\left[  \frac{1}{\left|  x+W_{T}\right|  ^{\left(
d-\alpha-1\right)  \left(  1+\delta\right)  p^{\prime}}}\right]  +\frac
{1}{\left|  x\right|  ^{\left(  d-\alpha-1\right)  \left(  1+\delta\right)
p^{\prime}}}\right\} \\
& \qquad\qquad +C\int_{0}^{T}E\left[  \frac{1}{\left|  x+W_{t}\right|  ^{\left(
d-\alpha+1\right)  \left(  1+\delta\right)  p^{\prime}}}\right]  dt.
\end{align*}
With the notation $p_{t}\left(  y\right)  =\frac{1}{\sqrt{\left(  2\pi\right)
^{d}t^{d}}}\exp\left(  -\frac{\left|  y\right|  ^{2}}{2t}\right)  $, and the
bound $\int_{0}^{T}p_{t}\left(  y\right)  dt\leq\frac{C_{d}\exp\left(
-\left|  y\right|  \right)  }{\left|  y\right|  ^{d-2}}$ we have%
\[
E\left[  \frac{1}{\left|  x+W_{T}\right|  ^{\left(  d-\alpha-1\right)  \left(
1+\delta\right)  p^{\prime}}}\right]  =\int_{\mathbb{R}^{d}}\frac{1}{\left|
x+y\right|  ^{\left(  d-\alpha-1\right)  \left(  1+\delta\right)  p^{\prime}}%
}p_{T}\left(  y\right)  dy
\]
and
\begin{align*}
&  \int_{0}^{T}E\left[  \frac{1}{\left|  x+W_{t}\right|  ^{\left(
d-\alpha+1\right)  \left(  1+\delta\right)  p^{\prime}}}\right]  dt\\
&\qquad\qquad  =\int_{\mathbb{R}^{d}}\frac{1}{\left|  x+y\right|  ^{\left(  d-\alpha
+1\right)  \left(  1+\delta\right)  p^{\prime}}}\left(  \int_{0}^{T}%
p_{t}\left(  y\right)  dt\right)  dy\\
&\qquad\qquad  \leq\int_{\mathbb{R}^{d}}\frac{1}{\left|  x+y\right|  ^{\left(
d-\alpha+1\right)  \left(  1+\delta\right)  p^{\prime}}}\frac{C_{d}\exp\left(
-\left|  y\right|  \right)  }{\left|  y\right|  ^{d-2}}dy.
\end{align*}

Thus, with a new constant $C>0$ depending on $\delta$ and the other
parameters,%
\[
I_{2}\leq C\left(  1+I_{2}^{\left(  1\right)  }+I_{2}^{\left(  2\right)
}+I_{2}^{\left(  3\right)  }\right)
\]
where%
\[
I_{2}^{\left(  1\right)  }:=\int_{\mathbb{R}^{d}}\int_{\mathbb{R}^{d}%
}e^{-\lambda_{\delta}\left|  x\right|  }\frac{1}{\left|  x+y\right|  ^{\left(
d-\alpha-1\right)  \left(  1+\delta\right)  p^{\prime}}}p_{T}\left(  y\right)
dxdy
\]%
\[
I_{2}^{\left(  2\right)  }:=\int_{\mathbb{R}^{d}}e^{-\lambda_{\delta}\left|
x\right|  }\frac{1}{\left|  x\right|  ^{\left(  d-\alpha-1\right)  \left(
1+\delta\right)  p^{\prime}}}dx
\]%
\[
I_{2}^{\left(  3\right)  }:=\int_{\mathbb{R}^{d}}\int_{\mathbb{R}^{d}%
}e^{-\lambda_{\delta}\left|  x\right|  }\frac{1}{\left|  x+y\right|  ^{\left(
d-\alpha+1\right)  \left(  1+\delta\right)  p^{\prime}}}\frac{C_{d}\exp\left(
-\left|  y\right|  \right)  }{\left|  y\right|  ^{d-2}}dxdy.
\]%
Choose $\delta>0$ such that $\left(  d-\alpha+1\right)  \left(  1+\delta
\right)  p^{\prime}<d$. Since $$\left(  d-\alpha-1\right)  \left(
1+\delta\right)  p^{\prime}< \left(  d-\alpha+1\right)  \left(
1+\delta\right)  p^{\prime}<d$$ the term $I_{2}^{\left(  2\right)  }$ is
finite. For $I_{2}^{\left(  1\right)  }$ and $I_{2}^{\left(  3\right)  }$ it
is sufficient to integrate first in $x$, bound the result uniformly in $y$,
then integrate in $y$; one proves that $I_{2}^{\left(  1\right)  }$ and
$I_{2}^{\left(  3\right)  }$ are finite. The proof is complete.

\section{The energy of a random vortex filament}

In~\cite{Ffil,FGub,Ascona} with the purpose of modeling turbulence
in 3d fluids, the authors introduce and study a model of random vortex
filaments based on Brownian motion. This model has been extended to
the fBm with $H>1/2$
by~\cite{Nua} and~\cite{FlandoliMinnelli}. Here we
recall briefly the model, emphasize the relationship of the vortex
energy with the pathwise regularity of the current associated with the
vortex \emph{core} and obtain new conditions for the integrability of the
vortex energy for the case $H \in (1/4,1/2)$.

\medskip For simplicity we consider only a single vortex since
extension to a linear superposition of different vortexes is
straightforward (and even to a random field of Poissonian vortexes,
see for example~\cite{FGubStat}). Let $(X_t)_{t\in[0,T]}$, $T>0$ be a 3d
fBm with Hurst parameter $H \in (1/4, 1)$ and consider the associated
vector current, formally written as
$$
\xi_{0}(x) = \int_0^T \delta(x-X_t) dX_t
$$
where the integral is a symmetric (Stratonovich) integral. This
object should be understood according to theorem~\ref{th:l2result} 
that is as a random distribution in the Sobolev space $H^{-\alpha}_{p}(\RR^d)$
of sufficiently large negative order. 
The vorticity field is then built by superposing translates of this
core weighted according to a compactly supported signed measure $\rho$
with finite mass which determines the intensity of vorticity. 
For more details about those considerations, the reader can consult 
\cite{Ffil}.
Then we end up with
$$
\xi(x) = \int_{\R^3}  \xi_0(x-y) \rho(dy)
$$
which is again a random distribution. 

The velocity field $u$ is generated from $\xi$ according to the Biot-Savart
relation
\begin{equation}
u(x)= \int_{\R^3} \mathcal{K}(x-y) \wedge \xi(y) dy
= \int_{\R^3} \mathcal{K}*\rho (x-y) \wedge \xi_0(y) dy 
\end{equation}
where $\wedge$ is the vector product in $\R^3$ and the vector kernel
$\mathcal{K}(x)$ is defined as
$\mathcal{K}(x):=(4\pi)^{-1} x/|x|^3
$ and $\mathcal{K}*\rho$ denote the convolution
$
(\mathcal{K}*\rho)(x) = \int_{\R^3} \mathcal{K} (x-z) \rho(dz).
$
The kinetic energy of the fluid is then defined as the $L^2(\RR^3)$ norm of $u$:
\begin{equation}
  \label{eq:energy}
\mathcal{E} = \int_{\R^3} |u(x)|^2 dx = \|u\|^2.  
\end{equation}
It is then interesting to find conditions on $\rho$ such that the kinetic
energy of the fluid is finite. Abstractly we have $u = \Phi \xi_0$
where we introduced an operator $\Phi$  whose kernel is
$\mathcal{K}*\rho$ having Fourier transform
$$
\mathcal{F}(\mathcal{K}*\rho)(q) = \frac{iq}{|q|^2}
\widehat \rho(q)
$$
where we denoted $\widehat \rho$ the Fourier transform of the
measure $\rho$. 

From now on $L^2$ will stay for $L^2(\RR^3)$.
Since, by Corollary~\ref{cor:current-l2-reg}, $\xi_0$ belongs a.s. to the
space $H^{-\alpha}_{2}(\RR^d)$ for any $\alpha > \alpha_H = 1/(2H) +
1/2$  (since $d=3$),
the condition $u \in L^2$ a.s. can be satisfied if $\Phi :
H^{-\alpha}_{2}(\RR^3) \to L^2$ which in Fourier variables is
sufficient to
require that
$$
\|\Phi\|_{H^{-\alpha}_{2}\to L^2} = \|\Phi (1-\Delta)^{\alpha/2}\|_{L^2
\to L^2} = \text{ess}\sup_{q \in \R^3}  
\frac{|\widehat \rho(q)|}{|q|} (1+|q|^2)^{\alpha/2} < \infty.
$$ 
for some  $\alpha > \alpha_H.$

We  can now formulate the following result.

\begin{corollary}
The kinetic energy of the vortex filament $\xi$ built upon a 3d
fractional Brownian motion of Hurst index $H>1/4$ is a.s. finite and in
$L^1$ if the measure $\rho$ satisfies
\begin{equation}
  \label{eq:our-finiteness}
\mathrm{ess}\sup_{q} |\widehat \rho(q)| |q|^{-1}(1+|q|^2)^{\alpha/2} < \infty  
\end{equation}
for some $\alpha > \alpha_H $.
\end{corollary}

\begin{remark}
\begin{enumerate}
\item
Known conditions on $\rho$ which guarantee the integrability of the
energy are given in~\cite{Ffil,FGub} for the case of Brownian motion
and It\^o  Brownian processes, and in~\cite{Nua} for the case
of fractional Brownian motions with Hurst parameter $H> 1/2$. From
\cite{Nua} it can be deduced 
that a sufficient condition for the integrability of the energy is
\begin{equation}
  \label{eq:nua-finiteness}
\int_{\R^3} dq \frac{|\widehat \rho(q)|^2}{|q|^{4-1/H}} <\infty  
\end{equation}
or, written in a different but equivalent form,
$$
\int_{\R^3}  \int_{\R^3} \frac{\rho(dx) \rho(dy)}{|x-y|^{1/H-1}} < \infty.
$$
\item  Condition (\ref{eq:our-finiteness})
implies Condition   (\ref{eq:nua-finiteness}) when $H > 1/2$.
 
In fact the left-hand side of 
 Condition~(\ref{eq:our-finiteness}) 
gives 
\begin{eqnarray*}
\int_{\R^3} dq |\widehat \rho(q)|^2   |q|^{-2}(1+|q|^2)^{\alpha} & & 
|q|^{2}(1+|q|^2)^{-\alpha} 
 |q|^{1/H -4} = 
\int_{\R^3} dq |\widehat \rho(q)|^2   |q|^{1/H - 2}
(1+|q|^2)^{-\alpha} \\
&\le&
 A \int_0^{+\infty} r^{1/H}  (1+r^2)^{-\alpha} dr ,
\end{eqnarray*}
where $A$ is the finite quantity of (\ref{eq:our-finiteness}).
Clearly previous expression is bounded for $\alpha  > \alpha_H$.

Of course the converse is not true.

\item
A way of finding similar conditions  to~(\ref{eq:nua-finiteness})   is to
follow the steps of Sec.~\ref{sec:hilbert} for  the Hilbert
space regularity of the stochastic currents and rewrite (formally) the
kinetic energy as   
\begin{equation}
  \label{eq:energy2}
\mathcal{E}= \int_0^T \int_0^T \langle  dX_t,  g(X_t-X_s)  dX_s\rangle
\end{equation}
where $g$ is a vector kernel with the following Fourier transform 
\begin{equation*}
  \begin{split}
\widehat g(q) &= \frac{|\widehat \rho(q)|^2}{|q|^2} \Pi_q
  \end{split}
\end{equation*}
and $\Pi_q$ is the following matrix
$$
(\Pi_q)_{\alpha \beta} = \delta_{\alpha \beta} - \frac{q_\alpha
q_\beta}{|q|^2}, \qquad \alpha,\beta = 1,\dots,3,
$$
which projects in directions orthogonal to $q$.
Formula~(\ref{eq:energy2}) can be understood, formally according to
Theorem~\ref{th:l2result}, as being the limit of the expectations   of
$\eps$-approximations 
$$ \mathcal{E}_\eps =
  \int_{0}^{T}\int_{0}^{T} g(  X_{t}-X_{s})
\left\langle D_{\varepsilon}X_{t},D_{\varepsilon}X_{s}\right\rangle
dtds 
$$
 To obtain conditions
for its finiteness in the spirit of eq.~(\ref{eq:nua-finiteness}), we
need to follow again the computations involved in the proof of
Theorem~\ref{th:reg-symm} and use a different strategy in bounding
some terms. Then we can prove the following:
\end{enumerate}
\end{remark}

\begin{theorem}
\label{th:symmetric-vortex}
Let $H>1/4$ and let $\rho: \RR^3 \to \RR$ be a function with Fourier transform $\widehat \rho$ satisfying
\begin{equation}
  \label{eq:vortex-cond}
\int_{\R^3} dq \frac{|\widehat \rho(q)|^2}{ |q|^{4-1/H}} <\infty  
\end{equation}
then the family of random fields $\{ u_\eps \}_{\eps\in(0,1)}$ defined as
$$
u_\eps(x) = \int_0^1 (\mathcal{K} *\rho)(x-X_t) \wedge D^0_\eps X_t dt, \qquad x
\in \R^3
$$
converges a.s. in $L^{2-\theta}(\Omega; L^2(\R^3;\R^3))$ for any
$\theta > 0$ to a random field $u \in L^2(\R^3;\R^3)$.
\end{theorem}

\begin{proof}
We will prove that $\sup_{\eps \in (0,1)} \expect \|u_\eps\|^2 <
\infty$ following the lines of the proof of Theorem~\ref{th:reg-symm}, then the
conclusion follows applying Theorem~\ref{th:l2result}. Let $\mathcal{E}_\eps = \|u_\eps\|^2$.

We start by treating the case $H > 1/2$.

Using Theorem \ref{th:wick} (Wick theorem)
 and independence of different coordinates we have
\begin{equation}
\label{eq:wick1ter}
  \begin{split}
\expect \mathcal{E}_\eps & =
\expect     
\int_{0}^{T}\int_{0}^{T} \sum_i g_{ii}\left(
X_{t}-X_{s}\right)  
\cov(D^0_\eps X^1_t, D^0_\eps X^1_s)
\,dt ds
\\ & \qquad -
\expect     
\int_{0}^{T}\int_{0}^{T} \sum_{ij} \nabla_i \nabla_j g_{ij}\left(
X_{t}-X_{s}\right)  
|\cov(D^0_\eps X^1_t,  X^1_{t}-X^1_{s})|^2
\,dt ds
\\  & =
\expect     
\int_{0}^{T}\int_{0}^{T}\text{Tr} g\left(
X_{t}-X_{s}\right)  
\cov(D^0_\eps X^1_t, D^0_\eps X^1_s)
\,dt ds
  \end{split}
\end{equation}
since a direct computation shows that
$
\sum_i \nabla_i g_{ik}(x) = 0
$.
Using the first bound in Lemma~\ref{lemma:bounds} we get

\begin{equation}
\label{eq:wick1bis}
  \begin{split}
\expect \mathcal{E}_\eps & \le \const
\expect     
\int_{0}^{T}\int_{0}^{t} \text{Tr} g\left(
X_{t}-X_{s}\right)  |t-s|^{2H-2}
\,ds dt
\\ & = \const \int_{\R^3} dq \text{Tr} \widehat g(q)
\int_{0}^{T}\int_{0}^{t} |t-s|^{2H-2} \expect  e^{-i \langle q,X_t-X_s\rangle}
\,ds dt
\\ & = \const \int_{\R^3} dq \text{Tr} \widehat g(q)
\int_{0}^{T}\int_{0}^{t} |t-s|^{2H-2} e^{-|q|^2 (t-s)^{2H}/2}
\,ds dt
\\ &\le \text{const}  \int_{\R^3}  dq |\widehat g(q)| \int_{0}^T dt 
\int_{0}^{\infty}
d\tau \tau^{2H-2}  e^{-|q|^2 \tau^{2H}/2}
\\ & =  \text{const}  \int_{\R^3}  dq |\widehat g(q)| |q|^{1/H-2} \int_{0}^T dt \int_{0}^{\infty}
dy y^{1-1/H}  e^{-y^2/2}
\\ & \le \text{const} T \int_{\R^3}  dq |\widehat g(q)| |q|^{1/H-2}
  \end{split}
\end{equation}
since
$$
\int_{0}^{\infty}
dy y^{1-1/H}  e^{-y^2/2} < \infty
$$
for $1-1/H > -1$, that is $H > 1/2$.  
Sufficient condition for uniform boundedness of $\expect \mathcal{E}_\eps$ is that
$$
\int_{\R^3} dq |\widehat g(q)| |q|^{1/H-2}  < \infty.
$$

Let us now consider the case $H \le 1/2$ and rewrite the approximated
energy as
$$
\mathcal{E}_\eps = \|u_\eps\|^2 = - \int_0^1 \int_0^1 \langle h(X_t-X_s)
D_\eps X_s,D_\eps X_t \rangle + \int_0^1 \int_0^1 \langle
g(0) D_\eps X_s,D_\eps X_t \rangle 
$$   
where $h(x) = g(0) - g(x)\ge 0$. Note that $g(0)$ is well defined using the
hypothesis of the theorem about the integrability of its Fourier
transform, moreover, as in Thm.~\ref{th:reg-symm} (in the $H<1/2$ part) we have the limit
 $$
 \int_0^1 \int_0^1 \langle
g(0) D_\eps X_s,D_\eps X_t \rangle \to   \langle (X_1-X_0)
g(0) ,(X_1-X_0) \rangle.
$$
So let us focus on the double integral with the kernel
$h$. Proceeding as in the $H > 1/2$ case we have
\begin{equation*}
  \begin{split}
 J & = - \int_0^T \int_0^T \langle h(X_t-X_s)
D_\eps X_s,D_\eps X_t \rangle
\\  &   \le - \const \int_0^T \int_0^T \int_{\R^3} dq  \text{Tr} \widehat h(q) \expect[ e^{i
\langle q,X_t-X_s \rangle} ] |t-s|^{2H-2 }\, dtds
\\  &   = - \const \int_0^T \int_0^T \int_{\R^3} dq  \text{Tr} \widehat h(q) e^{-|q|^2/2(t-s)^{2H}}
 |t-s|^{2H-2 }\, dtds
  \end{split}
\end{equation*}
but since $\widehat h(q) =  g(0) \delta(q) - \widehat g(q) $ we have
\begin{equation*}
  \begin{split}
& \int_{\R^3} dq  \text{Tr} h(q) e^{-|q|^2/2(t-s)^{2H}}
= \int_{\R^3} dq  \text{Tr} [g(0) \delta(q) - \widehat g(q)] e^{-|q|^2/2(t-s)^{2H}}
\\ & \qquad = \int_{\R^3} dq  \text{Tr} [g(0) \delta(q) - \widehat g(q)] [e^{-|q|^2/2(t-s)^{2H}}-1]
= - \int_{\R^3} dq  \text{Tr} \widehat g(q) [e^{-|q|^2/2(t-s)^{2H}}-1]    
  \end{split}
\end{equation*}
Then
\begin{equation*}
  \begin{split}
|J| & \le \const \int_0^T \int_0^T \int_{\R^3} dq  |\widehat g(q)| (1- e^{-|q|^2/2(t-s)^{2H}})
 |t-s|^{2H-2 }\, dtds    
\\&\le \text{const}  \int_{\R^3}  dq |\widehat g(q)| \int_{0}^T dt \int_{0}^{t}
ds (t-s)^{2H-2}  (1-e^{-2 |q|^2 (t-s)^{2H}})
\\ &\le \text{const}  \int_{\R^3}  dq |\widehat g(q)| \int_{0}^T dt \int_{0}^{\infty}
d\tau \tau^{2H-2}  (1-e^{-|q|^2 \tau^{2H}/4})
\\ & =  \text{const}  \int_{\R^3}  dq |\widehat g(q)| |q|^{1/H-2} \int_{0}^1 dt \int_{0}^{\infty}
dy y^{1-1/H}  (1-e^{-y^2/4}-1)
\\ & \le \text{const}  \int_{\R^3}  dq |\widehat g(q)| |q|^{1/H-2}
  \end{split}
\end{equation*}
where we made a change of variables $y = |q|\tau^H$ and we used the fact that
$$
\int_{0}^{\infty}
dy y^{1-1/H}  (1-e^{-y^2/4}) \le \int_{0}^{\infty}
dy y^{1-1/H}  \min(y^2,1) < \infty
$$
since $-3 < 1 - 1/H < -1$ when $1/4 < H < 1/2$.
So we obtain the uniform boundedness of $\expect \mathcal{E}_\eps$ when eq.~(\ref{eq:vortex-cond}) is satisfied. 

Analogously, the case $H=1/2$ does not pose any additional problem.
\end{proof}

\begin{remark}
Note that for $H \ge 1/2$ we recover condition~(\ref{eq:nua-finiteness}).
However while in~\cite{Ffil,Nua} only the existence and the
integrability properties of the energy are studied, here we have also 
 informations about convergence of $\eps$-approximations of
the velocity field generated by the random vortexes.  
\end{remark}

\appendix

\section{Some proofs and auxiliary results}
\label{sec:appA}

\begin{proof}[Proof of Lemma~\ref{lemma su K}]
Denote $\frac{1}{\Gamma\left(  \alpha\right)  \left(  4\pi\right)  ^{d/2}}$ by
$\gamma$, for shortness. Notice that for $x=0$ we have
\[
K_{\alpha}\left(  0\right)  =\gamma\int_{0}^{\infty}t^{\alpha-\frac{d}{2}%
}e^{-t}\frac{dt}{t}<\infty\text{ if and only if }\alpha>\frac{d}{2}.
\]
For $x\neq0$ we may use the change of variables $t=\left|  x\right|  ^{2}s$
and get
\begin{align*}
K_{\alpha}\left(  x\right)  =\left|  x\right|  ^{2\alpha-d}\rho\left(
x\right)  \\
\rho\left(  x\right)  :=\gamma\int_{0}^{\infty}s^{\alpha-\frac{d}{2}}%
e^{-\frac{1}{4s}-\left|  x\right|  ^{2}s}\frac{ds}{s}%
\end{align*}
where the integral converges for every value of the parameters, thanks to the
exponentials. For $0<\alpha<\frac{d}{2}$, we have
\[
c_{\alpha,d}e^{-2\left|  x\right|  ^{2}}\leq\gamma\int_{1}^{2}s^{\alpha
-\frac{d}{2}}e^{-\frac{1}{4s}-\left|  x\right|  ^{2}s}\frac{ds}{s}\leq
\rho\left(  x\right)
\]
for a positive constant $c_{\alpha,d}$. Moreover,
\begin{align*}
\rho\left(  x\right)    & \leq\gamma\int_{0}^{\frac{1}{\left|  x\right|  }%
}s^{\alpha-\frac{d}{2}}e^{-\frac{1}{4s}}\frac{ds}{s}+\gamma\int_{\frac
{1}{\left|  x\right|  }}^{\infty}s^{\alpha-\frac{d}{2}}e^{-\frac{\left|
x\right|  }{4}-\left|  x\right|  ^{2}s}\frac{ds}{s}\\
& \leq C_{\alpha,d}\int_{0}^{\frac{1}{\left|  x\right|  }}e^{-\frac{1}{8s}%
}\frac{ds}{s^{2}}+C_{\alpha,d}e^{-\frac{\left|  x\right|  }{4}}\int_{\frac
{1}{\left|  x\right|  }}^{\infty}s^{\alpha-\frac{d}{2}}\frac{ds}{s}\\
& \leq C_{\alpha,d}e^{-\frac{\left|  x\right|  }{8}}+C_{\alpha,d}%
e^{-\frac{\left|  x\right|  }{4}}\left|  x\right|  ^{\frac{d}{2}-\alpha}\leq
C_{\alpha,d}e^{-\frac{\left|  x\right|  }{8}}%
\end{align*}
for a positive constant $C_{\alpha,d}$ that we do not rename at every step.

If $\alpha>\frac{d}{2}$, we directly have from the original formula
\begin{align*}
c_{\alpha,d}^{\prime}e^{-\frac{\left|  x\right|  ^{2}}{8}}\leq\gamma\int
_{1}^{2}t^{\alpha-\frac{d}{2}}e^{-\frac{\left|  x\right|  ^{2}}{4t}-t}%
\frac{dt}{t}\leq K_{\alpha}\left(  x\right)  \\
\leq\gamma\int_{0}^{\left|  x\right|  }t^{\alpha-\frac{d}{2}}e^{-\frac{\left|
x\right|  ^{2}}{4t}}\frac{dt}{t}+\gamma e^{-\frac{\left|  x\right|  }{4}}%
\int_{\left|  x\right|  }^{\infty}t^{\alpha-\frac{d}{2}}e^{-t}\frac{dt}{t}\\
\leq C_{\alpha,d}^{\prime}e^{-\frac{\left|  x\right|  }{8}}%
\end{align*}
for a positive constants $c_{\alpha,d}^{\prime}$, $C_{\alpha,d}^{\prime}$. 
To estimate $K_{\alpha}(0)-K_{\alpha}(x)$ we write
\begin{equation*}
\begin{split}
K_{\alpha}(0)-K_{\alpha}(x) & = \gamma \int_{0}^{\infty}t^{\alpha-\frac{d}{2}}e^{-t} (1-e^{-\frac{\left|
x\right|  ^{2}}{4t}})\frac{dt}{t}
\\ & = \gamma |x|^{\alpha-\frac d2} \int_{0}^{\infty}s^{\alpha-\frac{d}{2}}e^{-s|x|^{2}} (1-e^{-\frac{1}{4s}})\frac{ds}{s}
\end{split}
\end{equation*}
and use the same arguments as above. In the case $\alpha = d/2$ we simply split the integral as above and by straightforward estimation we can prove that $K_\alpha(x) \le \const\ \log|x|$ for small $|x|$ and that $K_\alpha(x)$ decay exponentially for large $|x|$.
The
proof is complete.
\end{proof}

\medskip

\noindent
\begin{proof}[Proof of Lemma~\ref{lemma:deltaK}]
We observe that, since $K_{\alpha}$ is the kernel of the operator $(1-\Delta)^{-\alpha}$ we have the identity
$
K_{\alpha-1}(x) =  (1-\Delta) K_{\alpha}(x)
$ so that
$$
-\Delta K_{\alpha}(x) = K_{\alpha-1}(x) - K_{\alpha}(x) \qquad x \neq 0.
$$
Then $
|-\Delta K_{\alpha}(x)| \le| K_{\alpha-1}(x)| +| K_{\alpha}(x)|$ which gives the required bound using Lemma~\ref{lemma su K}. 
\end{proof}

\medskip

\noindent
\begin{proof}[Proof of Lemma~\ref{lemma:bounds}]
We start with the first estimate in 1. 
A direct computation shows
\begin{equation*}
\cov(D^0_\eps X^i_t, D^0_\eps X^i_s) = \frac{|t-s|^{2H-2}}{2} \Phi\left(\frac{2 \eps}{t-s}\right)
\end{equation*}
where
$$
\Phi(x) = \frac{|1+x|^{2H}+|1-x|^{2H}-2}{x^2}.
$$

The function $\Phi$ is continuous in $(0,\infty)$, $\lim_{x\to
0}\Phi(x) = 2H-1$ so, when $|t-s| \le 2\eps$ we have
$$
\frac{|t-s|^{2H-2}}{2} \Phi\left(\frac{2 \eps}{t-s}\right) \le \const\, |t-s|^{2H-2}.
$$
Moreover $\lim_{x\to \pm\infty} |x|^{2-2H} \Phi(x) = 2$, so when $|t-s|> 2\eps$ there exists a constant not depending on $\eps$ such that
$$
\frac{|t-s|^{2H-2}}{2} \Phi\left(\frac{2 \eps}{t-s}\right) \le \const\ \eps^{2H-2} \le \const\ |t-s|^{2H-2}
$$
which proves the first claim.

We discuss now the second estimate in 1. and point 2.
A direct computation gives
\begin{eqnarray*}
\cov(D^0_\eps X^i_t,  X^i_{t}-X^i_{s}) &=&  -
\cov(D^0_\eps X^i_s,  X^i_{t}-X^i_{s}) =
 \frac{1}{2\eps}\left(|t-s+\eps|^{2H}-|t-s-\eps|^{2H}\right) \\
&=& |t-s|^{2H - 1} \psi \left(\frac{\eps}{t-s}\right)
\end{eqnarray*}
where 
$ \psi(x) = \frac{|1+x|^{2H} - |1-x|^{2H}}{2x}$.
It is easy to show that $\psi$ is continuous and $\psi(0+) = 2H$, moreover $|x|^{2-2H} \psi(x) \to 2H$ when $x \to \pm\infty$.
This allows to conclude the proof.
\end{proof}

\medskip

\noindent
\begin{proof}[Proof of Lemma~\ref{lemma:bounds-forw}]
The first estimate in 1.
 is very similar to the previous lemma. 
A direct computation shows
\begin{equation} \label{EPhi-forw}
\cov(D^-_\eps X^i_t, D^-_\eps X^i_s) = \frac{|t-s|^{2H-2}}{2} \Phi\left(\frac{ \eps}{t-s}\right)
\end{equation}
where
$\Phi$
is the same as previously.
As for the other points
a direct computation gives
\begin{eqnarray*} 
\cov(D^-_\eps X^i_t,  X^i_{t}-X^i_{s}) &=&  -
\cov(D^-_\eps X^i_s,  X^i_{t}-X^i_{s}) =
 \frac{1}{2\eps}\left(|t-s+\eps|^{2H}-|t-s-\eps|^{2H} - \eps^{2H}\right) \\
&=& (t-s)^{2H - 1} \tilde \psi \left(\frac{\eps}{t-s} \right)
\end{eqnarray*}
where 
$ \tilde \psi(x) = \frac{x^{2H} +1 - (1-x)^{2H}}{2x}$.
Then, it is not difficult, arguing as in lemma~\ref{lemma:bounds} to
conclude.
In particular one can evaluate the limit in (\ref{EPhi-forw}).

\end{proof}



\begin{thebibliography}{99}

\bibitem{Alos}
 E. Alos, J. A. Leon, D. Nualart,  Stochastic Stratonovich calculus
 for fractional Brownian motion with Hurst parameter less than
 $1/2$.  \it{Taiwanese J. Math.}  {\bf 5}  (2001),  no. 3, 609--632.  


\bibitem{ELY}K.D. Elworthy, Xue-Mei Li, M. Yor, The importance of strictly
local martingales;\ applications to radial Ornstein-Uhlenbeck processes,
\textit{Probab. Theory Relat. Fields} {\bf 115} (1999), 325-355.

\bibitem{em}  P. Embrechts, P., M. Maejima,  Selfsimilar processes.
 \textit{Princeton Series in Applied Mathematics}.  NJ, 2002. 


\bibitem{Ffil}  F. Flandoli, On a probabilistic description of small scale
structures in 3D fluids, \textit{Annales Inst. Henri
Poincar\'{e}, Probab. \& Stat.} \textbf{38} (2002), 207-228.

\bibitem{FGub}  F. Flandoli, M. Gubinelli, The Gibbs ensemble of a vortex
filament, \textit{Probab. Theory Rel. Fields} \textbf{122} (2002), 317-340.

\bibitem{FGubStat}  F. Flandoli, M. Gubinelli, Statistics of a vortex filament model, 
\textit{Electr. J. Prob.} \textbf{10} (2005), no. 25, 865--900. 



\bibitem{Ascona}  F. Flandoli, M. Gubinelli, Random Currents and
  Probabilistic Models of Vortex Filaments, \textit{
  Seminar on Stochastic Analysis, Random Fields and
  Applications IV, Ascona 2002}, Prog. In Prob. 58, Birk\"auser
  Verlag, Basel 2004. 


\bibitem{FlandoliMinnelli} F. Flandoli, I. Minelli,
 Probabilistic models of vortex filaments.
\textit{ Czechoslovak Math. J.}  \textbf{51} (126)  (2001),  no. 4, 713--731.

\bibitem {FGGT} F. Flandoli, M. Giaquinta, M. Gubinelli and V. M. Tortorelli,
Stochastic currents, \textit{Stoch. Proc. Appl.}
 \textbf{115}  (2005),  no. 9, 1583--1601.

\bibitem {gall} J.-F., Le Gall, Sur le temps local d'intersection du mouvement brownien plan et la m\'ethode
 de renormalisation de Varadhan.  S\'eminaire de probabilit\'es, XIX, 1983/84, 
 314--331, Lecture Notes in Math., 1123, Springer, Berlin, 1985.

\bibitem{Grad} M. Gradinaru, I. Nourdin, F. Russo, P. Vallois,
$m$-order integrals and generalized Ito's formula:
 the case of a fractional  Brownian motion with any Hurst index,
 \textit{Ann. Inst. H. Poincar\'e Probab. Statist.}  \textbf{41}  (2005),  no. 4, 781--806.

\bibitem{GradNo} M. Gradinaru, I. Nourdin, 
 Approximation at first and second order of $m$-order integrals of the fractional Brownian 
motion and of certain semimartingales, \textit{ Electron. J. Probab.}
\textbf{ 8}  (2003), no. 18, 26 pp.

\bibitem {Gubi} M. Gubinelli, Controlling rough paths, \textit{ J. Funct. Anal.}
\textbf{216} (2004), \textit{no. }1, 86--140.

\bibitem{qianlyons}   T. J. Lyons, Z. Qian, 
System control and rough paths,
\textit{ Oxford University Press}, Oxford, 2002.

\bibitem{Lyons}  T. J. Lyons, Differential equations driven by rough
signals, \textit{Revista Math. Iberoamericana}, \textbf{14} (1998),
\textit{no.} 2, 215--310.

\bibitem{Nua}
D. Nualart, C. Rovira, and S. Tindel,
 Probabilistic models for vortex filaments based on fractional
  Brownian motion,
\textit {RACSAM Rev. R. Acad. Cienc. Exactas F\'\i s. Nat. Ser. A Mat.},
  {\bf 95} (2001), 213--218.

\bibitem{Pa}  A. Pazy, Semigroups of Linear Operators and
Applications to Partial Differential Equations,\textit{
 Springer-Verlag, New York},
1983.

\bibitem{ry} D. Revuz and M. Yor, Continuous martingales and Brownian motion.
\textit{ Third edition.
{\bf 293}. Springer-Verlag, Berlin}, 1999.

\bibitem{RV}
F. Russo, P. Vallois,
 Stochastic calculus with respect to continuous finite quadratic variation processes.
\textit{ Stochastics Stochastics Rep.}  {\bf 70}  (2000),  no. 1-2, 1--40.

\bibitem{RVSem}
 F. Russo, P. Vallois,
Elements of stochastic calculus via regularization.
Preprint LAGA, 2004-28.
http://front.math.ucdavis.edu/math.PR/0603224.
To appear: S\'eminaire de Probabilit\'es XL,
Eds. C. Donati-Martin, M. Emery, A. Rouault, C. Stricker.


\bibitem{Triebel} H. Triebel, Interpolation Theory, Function Spaces,
 Differential
Operators, \textit{North Holland, Amsterdam} 1978.

\end{thebibliography}
\end{document}